\documentclass[multi,makeidx]{cambridge7A}
  \usepackage[numbers]{natbib}
  \usepackage{amsmath,amsthm,amsfonts}
  \usepackage{graphicx}
  \usepackage{makeidx}
  \usepackage[usenames,dvipsnames]{color}

\numberwithin{equation}{chapter}

  \usepackage{comment}

  \usepackage{caption}
\usepackage[labelfont=rm]{subcaption}

\usepackage[pagebackref=true, colorlinks]{hyperref}

\hypersetup{pdffitwindow=true,linkcolor=blue,citecolor=blue,urlcolor=blue,menucolor=blue}

\newcommand{\bo}{{\mathbf 1}}
\newcommand{\gd}{{\delta}}
\newcommand{\gb}{{\beta}}
\newcommand{\figref}[1]{\hyperref[#1]{Figure \ref{#1}}}
\newcommand{\figsref}[1]{\hyperref[#1]{Figures \ref{#1}}}
\newcommand{\secref}[1]{\hyperref[#1]{\S \ref{#1}}}
\newcommand{\thmref}[1]{\hyperref[#1]{Theorem \ref{#1}}}
\newcommand{\conjref}[1]{\hyperref[#1]{Conjecture \ref{#1}}}

\newcommand\pf{\begin{proof}}
\newcommand\epf{\end{proof}}

\newcommand\be{\begin{equation}}
\newcommand\ee{\end{equation}}
\newcommand{\<}{\left\langle}
\renewcommand{\>}{\right\rangle}

\DeclareMathOperator{\mult}{mult}   
\DeclareMathOperator{\vol}{vol} 
\DeclareMathOperator{\discr}{discr} 
\DeclareMathOperator{\sgn}{sgn} 
\DeclareMathOperator{\tr}{tr} 
\DeclareMathOperator{\SL}{SL} 
\DeclareMathOperator{\GL}{GL} 
\DeclareMathOperator{\PSL}{PSL} 
 
\newcommand{\fC}{\mathfrak{C}} 
\newcommand{\cF}{\mathcal{F}} 
 
\newcommand{\cX}{\mathcal{X}} 
\newcommand{\cY}{\mathcal{Y}} 
\newcommand{\bk}{\backslash}
\newcommand{\Z}{\mathbb{Z}}
\newcommand{\R}{\mathbb{R}}
\newcommand{\C}{\mathbb{C}}
\newcommand{\N}{\mathbb{N}}
\newcommand{\Q}{\mathbb{Q}}
\newcommand{\bH}{\mathbb{H}}

\newcommand{\G}{\Gamma}      
\newcommand{\g}{\gamma}    
\newcommand{\gl}{\lambda}     
\newcommand{\ga}{\alpha}    
\newcommand{\gs}{\sigma}      
\newcommand{\gt}{\theta}      
\newcommand{\gz}{\zeta}      
      
\newcommand{\gep}{\epsilon}  
\newcommand{\vep}{\varepsilon} 
\newcommand{\mattwo}[4]
{\left(\begin{array}{cc}
                        #1  & #2   \\
                        #3 &  #4
                          \end{array}\right) }
\newcommand{\sC}{\mathscr{C}}
\newcommand{\dd}{\partial}
\newcommand{\twocase}[5]{#1 \begin{cases} #2 & \text{#3}\\ #4
&\text{#5} \end{cases}   }
\newcommand{\mattwos}[4]
{\bigl( \begin{smallmatrix}
                        #1  & #2   \\
                        #3 &  #4
\end{smallmatrix} \bigr)
}

\usepackage{natbib}
\usepackage{mathrsfs}
\input xy
\usepackage{amsbsy,graphicx}

\usepackage{geometry}
\usepackage{amssymb}
\usepackage{psfrag}
\usepackage{bbm}
\usepackage[all]{xy}

\makeindex

  \newtheorem{theorem}{Theorem}[section]

  \newtheorem{conjecture}[theorem]{Conjecture}

  \theoremstyle{definition}

  \theoremstyle{remark}

  {\bfseries}{\rmfamily}

\begin{document}

\title[Proceedings of the Durham Easter School 2014]{Dynamics and
  Analytic Number Theory}

\author{D. Badziahin, A. Gorodnik, N. Peyerimhoff}

\setcounter{tocdepth}{1}




\author[A. Kontorovich]{Alex Kontorovich (Rutgers University)\footnotemark}

\chapter{Applications of Thin Orbits}

\footnotetext[1]{110 Frelinghuysen Rd\\ Piscataway, NJ 08854\\
email: \texttt{alex.kontorovich@rutgers.edu}\\ 
url: \texttt{math.rutgers.edu/$\sim$alexk}\\
The author is partially supported by an NSF CAREER grant  DMS-1455705, an NSF FRG grant DMS-1463940, an Alfred P. Sloan Research Fellowship, and a BSF grant.}

\contributor{Alex Kontorovich
  \affiliation{Rutgers University}}

\textbf{Abstract} \quad 
This text is
based on a series of  three 
expository
lectures
  on a variety of topics 
  related to ``thin orbits,''
as  delivered 
at Durham University's Easter School on ``Dynamics and Analytic Number Theory'' in April 2014. 
 The first 
 lecture
 reviews
 closed geodesics on the modular surface and 
  the reduction theory of binary quadratic forms before discussing Duke's equidistribution theorem (for indefinite classes). The second lecture exposits three  quite different but
  (it turns out)
   not unrelated problems, %
     due to  Einsiedler-Lindenstrauss-Michel-Venkatesh,  McMullen, and Zaremba. The third lecture reformulates these in terms of the aforementioned thin orbits, and shows how all three would follow from a single  ``Local-Global'' Conjecture of  Bourgain and the author. We also describe some partial progress on the conjecture, which has lead to some 
results 
on
the original problems.


\section{Lecture 1: Closed Geodesics,
Binary Quadratic Forms, and Duke's Theorem\index{Duke's Theorem}}

This first  lecture has three parts.
In \secref{sec:closedGeo}, we review
the geodesic flow on the hyperbolic plane to
study closed geodesics on the modular surface.
Then \secref{sec:Gauss} discusses
  Gauss's reduction theory of
  binary quadratic forms.
  Finally, in \secref{sec:Duke}, we combine
the previous two discussions to 
connect  indefinite classes
  to
  closed geodesics, and
state
 Duke's equidistribution theorem. 

\subsection{Closed Geodesics}\label{sec:closedGeo}

 \begin{figure}
        \begin{subfigure}[t]{0.45\textwidth}
                \centering
		\includegraphics[width=\textwidth]{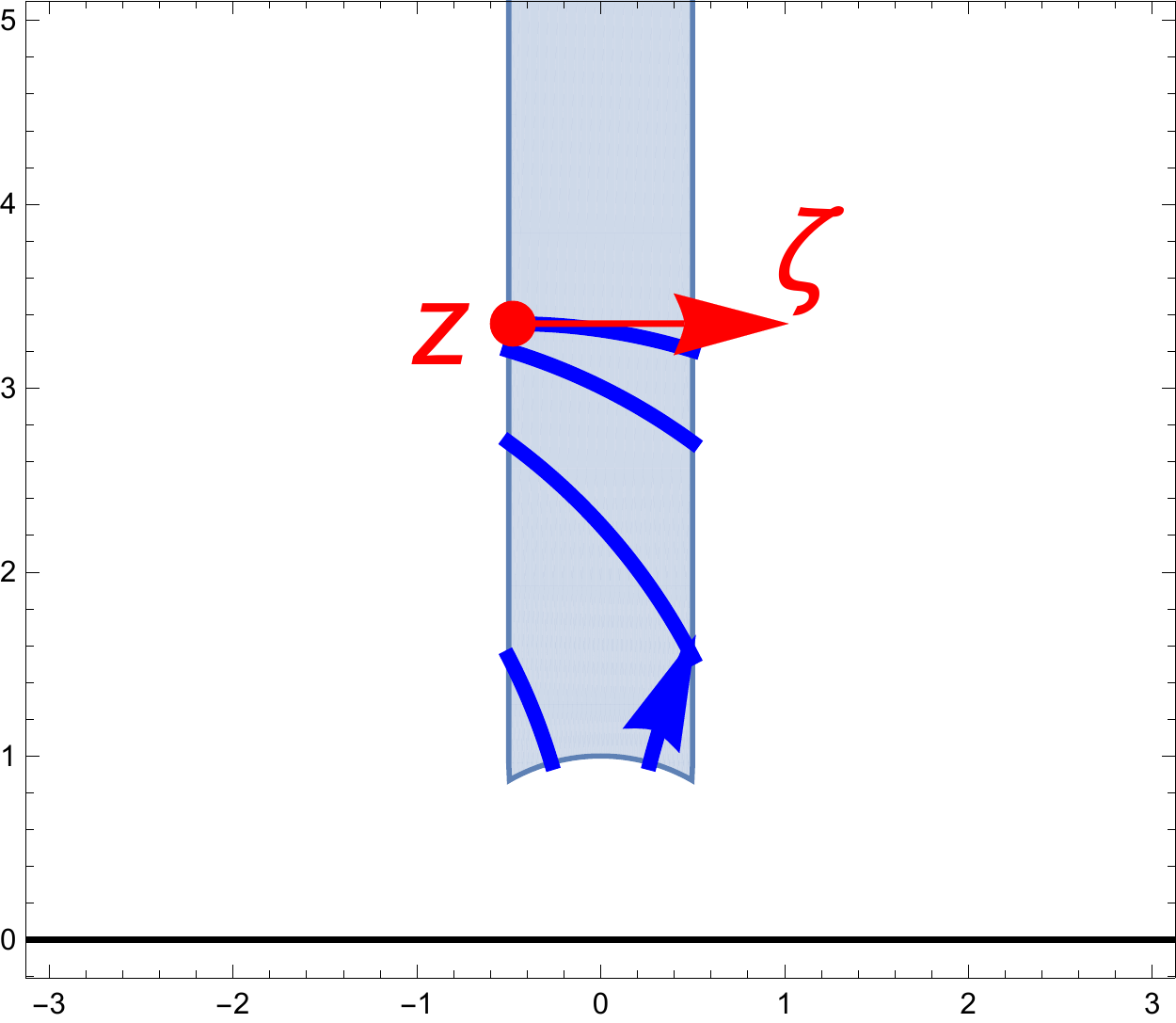}
                \caption{As a broken ray}
                \label{AKfig:1a}
        \end{subfigure}%
\qquad
        \begin{subfigure}[t]{0.45\textwidth}
                \centering
		\includegraphics[width=\textwidth]{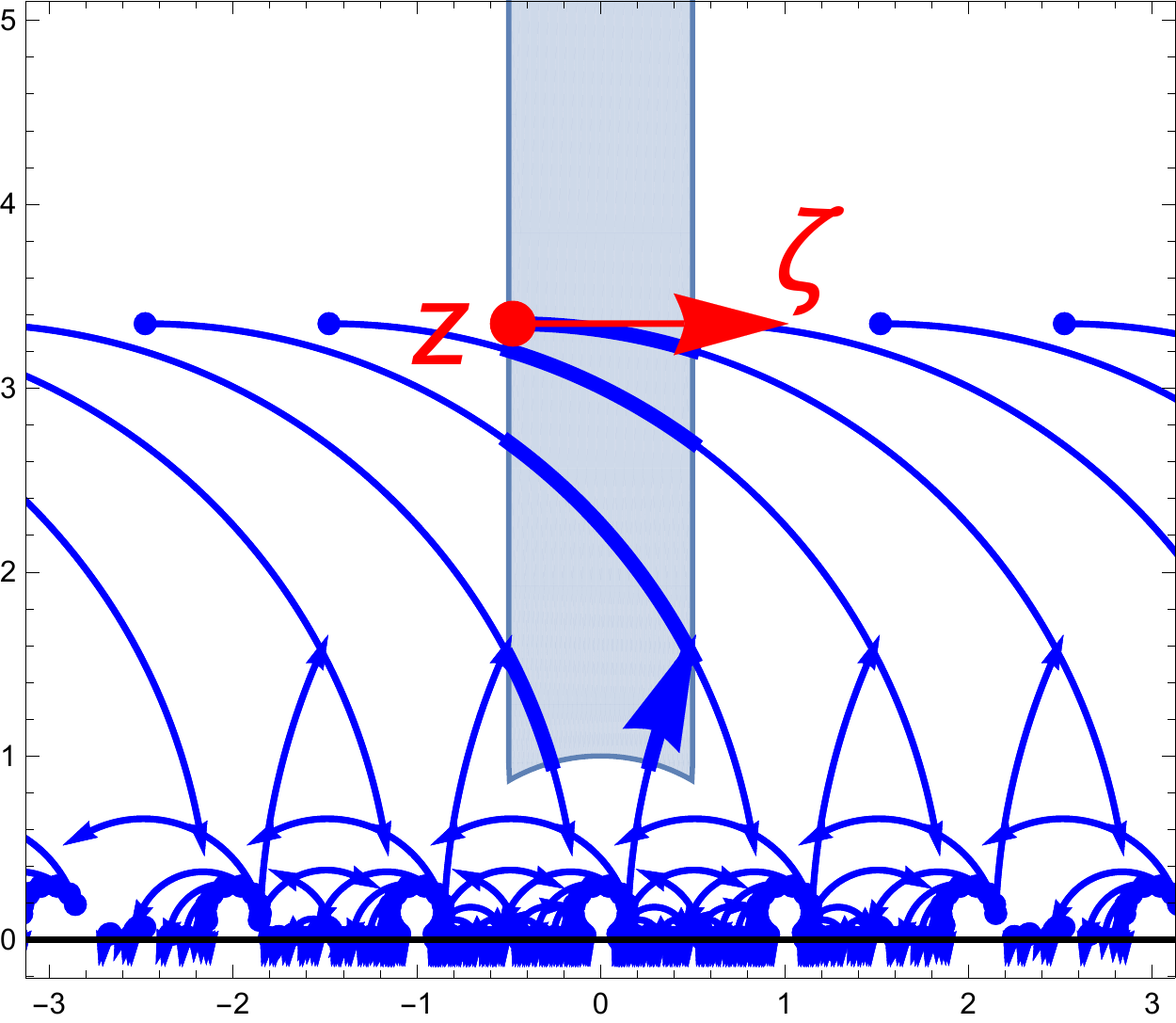}
                \caption{As the collection of $\PSL_{2}(\Z)$-translates} 
                \label{AKfig:1b}
        \end{subfigure}
\caption{The geodesic flow on the modular surface} 
\end{figure}

Let $\bH=\{x+iy:y>0\}$ denote the Poincar\'e (or perhaps more precisely, Beltrami) upper half plane,
let $T\bH$ be its tangent bundle, and
 for $(z,\gz)\in T\bH$, equip the tangent bundle 
with Riemannian metric $\|\gz\|_{z}:=|\gz|/\Im z$.
Here $z\in\bH$ is the ``position'' and $\gz\in T_{z}\bH\cong\C$ is the ``direction'' vector.
Let
  $T^{1}\bH$ be the {\it unit} tangent bundle
  of all $(z,\gz)\in T\bH$ having $\|\gz\|_{z}=1$.
The
fractional linear action\index{fractional linear action} of
   the group
 $G=\PSL_{2}(\R)=\SL_{2}(\R)/\{\pm I\}$ 
on $\bH$ induces the following action on $T^{1}\bH$:
\be\label{AKeq:Gact}
G\ni \mattwo abcd \  : \ (z,\gz)\in T^{1}\bH \ \mapsto \ \left({az+b\over cz+d},{\gz\over (cz+d)^{2}}\right)
,
\ee
with invariant measure
\be\label{AKeq:Haar}
d\mu={dx\, dy\, d\gt\over y^{2}},
\ee
in coordinates $(x+iy,\gz)$, where $\arg\gz=\gt$.

{\bf Exercise:} This is indeed an   action,
which is moreover free, transitive, and invariant for the measure in \eqref{AKeq:Haar}. 

The geodesics on $\bH$ are vertical half-lines and semi-circles orthogonal to the real line. 
Given $(z,\gz)\in T^{1}\bH$, the time-$t$ geodesic flow\index{geodesic flow} moves $z$ along the geodesic determined by $\gz$ to the point at distance $t$ from $z$.
The 
{\it visual point} from some $(z,\gz)\in T^{1}\bH$ is
the point on the boundary $\dd\bH\cong
\R\cup\{\infty\}
$ that one obtains 
by following the geodesic flow for
infinite time.

People studied such flows on  manifolds purely geometrically for 
some
time before Gelfand championed the injection of algebraic and representation-theoretic ideas. With Fomin \cite{GelfandFomin1951}, he discovered that under the  identification 
\be\label{AKeq:gIden}
G\cong T^{1}\bH, 
\qquad
g\leftrightsquigarrow g(i,\uparrow),
\ee 
the geodesic flow on $T^{1}\bH$ corresponds in $G$ to right multiplication by the diagonal subgroup 
$$
A\ =\ \left\{a_{t}:=\mattwo {e^{t/2}}{}{}{e^{-t/2}}\right\},
$$  
see, e.g., \cite{BekkaMayer2000, EinsiedlerWard2011}.
\\

To study (primitive, oriented) {\it closed} geodesics, we move to the
 modular surface\index{modular surface}, defined as the quotient $\G\bk\bH$, with $\G=\PSL_{2}(\Z)$. Its unit tangent bundle $\cX:=T^{1}(\G\bk\bH)$ is, as above, identified with $\G\bk G$, and the geodesic flow 
 again
 corresponds to right multiplication by $a_{t}$.
 It is useful to think of this flow in two equivalent ways: (i) as a broken ray in a fundamental domain for $\cX$ which is sent back inside when it tries to 
 exit
 (see \figref{AKfig:1a}), or (ii) as a whole collection of $\G$-translates of a single geodesic ray  in the universal cover, $\bH$, as in \figref{AKfig:1b}.
Thus we will sometimes write $g\in \G\bk G$ for the first notion, and $\G g$ for the second. 

To obtain a 
closed geodesic\index{closed geodesic} on the modular surface, we 
start at some point $\G g$ and come back to the same exact point (including the tangent vector) after a (least) time $\ell>0$ (here $\ell$ is for {\it length}). 
That is,
$$
\G g a_{\ell}\ = \ \G g,\qquad\text{ or equivalently, }\qquad
g a_{\ell}\ = \ M g,
$$
for some matrix\footnote{Technically, as $\G=\SL_{2}(\Z)/\{\pm I\}$, we should be using cosets $\pm M$ here. We will abuse notation
and treat elements of $\G$ as matrices, with the convention that their trace is positive.}
 $M\in\G$.
Then
\be\label{AKeq:Mconj}
M \ =\ g a_{\ell}g^{-1},
\ee
so $M$ has eigenvalues $e^{\pm\ell/2}$, and $g$ is a matrix of eigenvectors. Note that $M$ is hyperbolic with trace
$$
\tr M \ = \ 2\cosh(\ell/2),
$$
and the expanding eigenvalue, $\gl$, say, is given by
\be\label{AKeq:eval}
\gl \ = \ 
e^{\ell/2} \ = \ (\tr M+\sqrt{\tr M^{2}-4})/2.
\ee
Actually, since $\G g$ is only determined up to left-$\G$ action, the matrix  $M$ is only determined up to $\G$-conjugation (which of course leaves invariant its trace). In this way, primitive closed geodesics correspond to 
primitive (meaning not of the form $[M_{0}^{n}]$ for some $M_{0}\in\G$, $n\ge2$)
hyperbolic conjugacy classes $[M]$ in $\G$. 

We note already from \eqref{AKeq:eval} that the lengths of closed geodesics are far from arbitrary; since $M\in\G$, its eigenvalues are quadratic irrationals.
It will also be useful later to note the visual point of $g$.
Writing $M$ as  $M=\mattwos abcd$, then
if $c>0$,
the matrix $g$ of eigenvectors can be given by 

\noindent
{\bf Exercise}:
\be\label{AKeq:gIs}
g\ = \ \frac1{(c^{2}(\tr^{2}M-4))^{1/4}}\mattwo {\gl-d}{1/\gl-d}cc
.
\ee
The scaling factor is to ensure $g$ has determinant $1$.
If $c<0$, negate the first column in \eqref{AKeq:gIs}.
Any other choice of $g$ is obtained by rescaling the first column by a factor $\gs e^{t/2}$ and the second by $\gs e^{-t/2}$, $\gs\in\{\pm1\}$; this of course corresponds to the right action by $a_{t}$ in $\PSL_{2}(\R)$.
The visual 
point $\ga$ from $g$ is determined by computing
\be\label{AKeq:gaIs}
\ga\ = \
\lim_{t\to\infty}ga_{t}\cdot i \  = \ {\gl-d\over c} = {a-d+\sqrt{\tr^{2}M-4}\over 2c}.
\ee
Note again that this is a quadratic irrational, and its Galois conjugate $\overline\ga$ is the visual point of the backwards geodesic flow.
Note also that $\ga$ is independent of the choice of $g$ above.
Finally, we record here that  the fractional linear action of $M$ on $\R$ fixes $\ga$; indeed, starting from $g a_{\ell}=Mg$, multiply both side on the right by $a_{t}$, have that matrix act on the left by $i$, and take the limit as $t\to\infty$:
\be\label{AKeq:MgaEqga}
\ga \ = \ \lim_{t\to\infty}ga_{\ell}a_{t}\cdot i \ = \ \lim_{t\to\infty}Mg a_{t}\cdot i  \ = \ M\ga.
\ee

\begin{figure}
		\includegraphics[width=.75\textwidth]{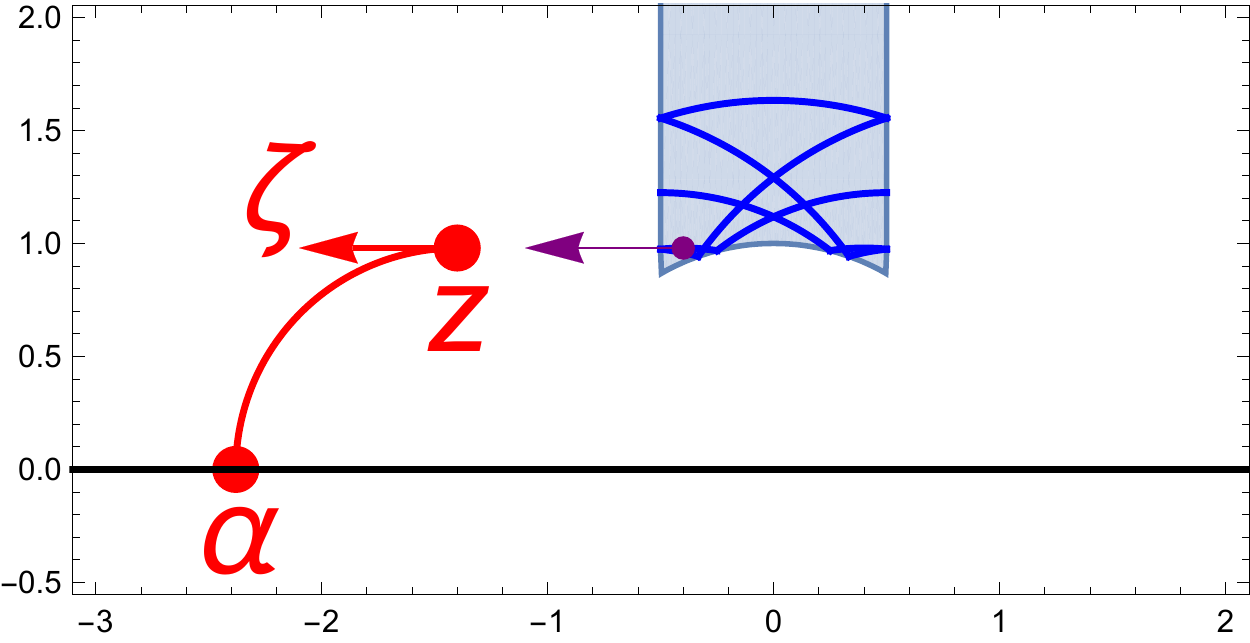}
\caption{The closed geodesic corresponding to $M$ in \eqref{AKeq:Mis}} 
\label{AKfig:geodM}
\end{figure}

To see an explicit example, let us construct the geodesic corresponding to the hyperbolic matrix 
\be\label{AKeq:Mis}
M=\mattwo {12}5{-5}{-2}.
\ee 
From \eqref{AKeq:eval}, \eqref{AKeq:gIs}, and \eqref{AKeq:gaIs}, we compute
\be\label{AKeq:gaIsEg}
\gl \ = \ 5+2 \sqrt{6},
\quad
g\ =\ 
\frac1{ \sqrt[4]{2400}}
\left(
\begin{array}{cc}
 -7-2 \sqrt{6} & 7-2 \sqrt{6} \\
 5 & -5 \\
\end{array}
\right)
,
\quad
\ga \ = \ \frac{-7-2 \sqrt{6}}{5}
.
\ee
Using the identification \eqref{AKeq:gIden} and action \eqref{AKeq:Gact}, the point $g\in G$ corresponds to the point $(z,\gz)\in T^{1}\bH$ where
\be\label{AKeq:zgzEg}
z \ = \
-\frac{7}{5}+\frac{2 i \sqrt{6}}{5}
,
\quad
\gz \ = \
-\frac{2 \sqrt{6}}{5}.
\ee
The points $z,\gz,$ and $\ga$ are shown in \figref{AKfig:geodM}, as well as their images in the standard fundamental domain 
$\cF$
for $\G$. The resulting closed geodesic, 
also
shown in
 $\cF$, has length $\ell=2\log \gl\approx4.58$.
Had we started with $M^{2}$ instead of $M$, we would have obtained the same $g, z, \gz$, and $\ga$ but the length would double, corresponding to looping around the geodesic twice (hence our restriction to primitive geodesics and conjugacy classes).
Replacing $M$ by some conjugate, $\g M\g^{-1}$, with $\g\in\G$, results in replacing $g, z$, and $\ga$ by $\g g$, $\g z$, and $\g \ga$, respectively (of course the geodesic remains unchanged).

Next we wish to discuss the cutting sequence\index{cutting sequence} of the geodesic flow. 
Recall that $\G$ is generated by two elements $T=\mattwos1101$ and $S=\mattwos01{-1}0,$ and that its standard fundamental domain $\cF$ is the intersection 
of the domains 
$$
\cF_{T}
\ := \ \{
\Re z>-1/2\},\
\cF_{T^{-1}}
\ := \ \{
\Re z<1/2\},
\text{ and }
\cF_{S} \ : = \
\{
|z|>1\}. 
$$
As we follow the geodesic flow from $\cF$, thought of as a subset of the universal cover, $\bH$, we pass through one of the boundary walls, leaving one of the domains $\cF_{L}$, $L\in\{T,T^{-1},S\}$; here $L$ is the ``letter'' we must apply to return the flow to $\cF$. Given a starting point $(z,\gz)\in T^{1}(\G\bk\bH)$, its cutting sequence is this sequence of letters $L$.

 \begin{figure}
        \begin{subfigure}[t]{0.3\textwidth}
                \centering
		\includegraphics[width=\textwidth]{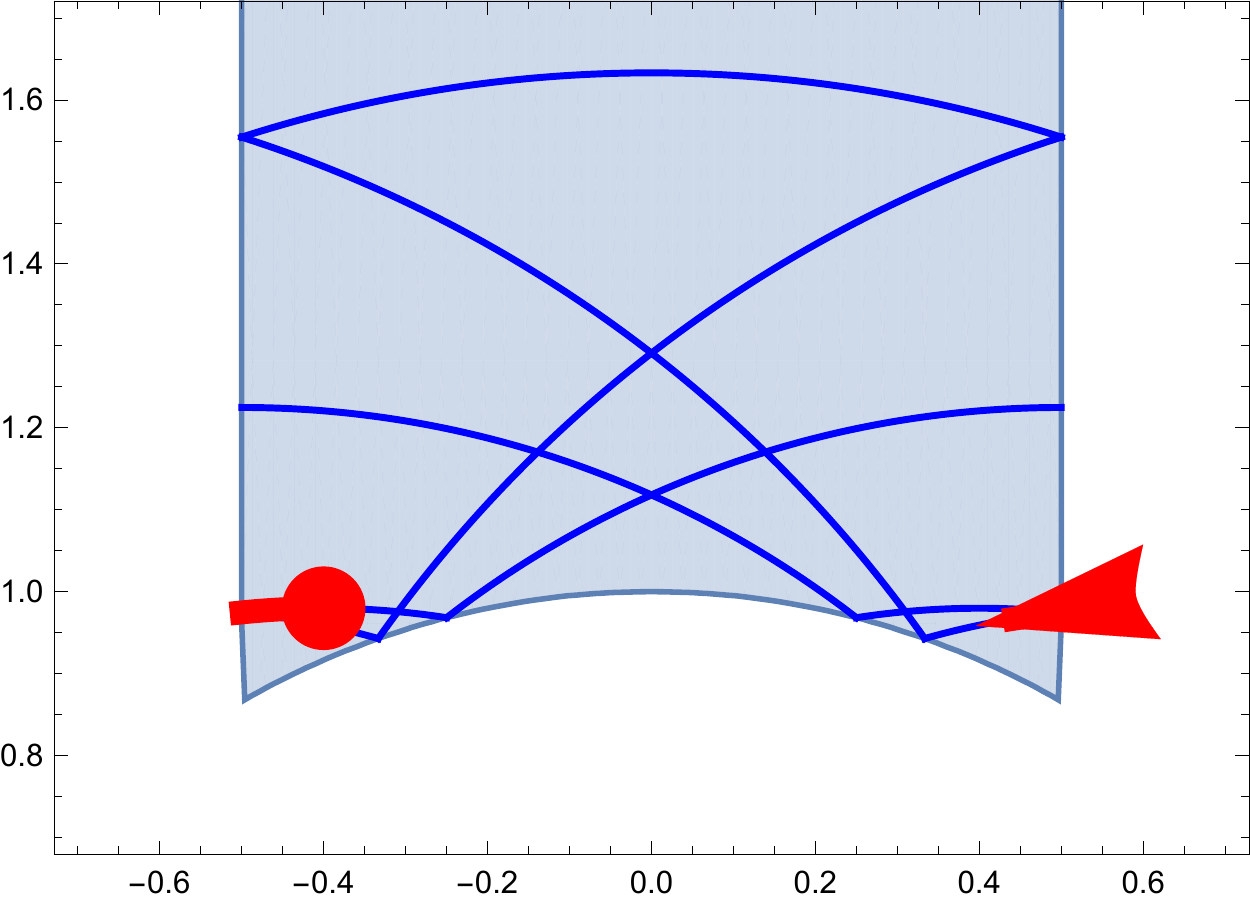}
                
                {$L_{1}=T$}
        \end{subfigure}%
\quad
        \begin{subfigure}[t]{0.3\textwidth}
                \centering
		\includegraphics[width=\textwidth]{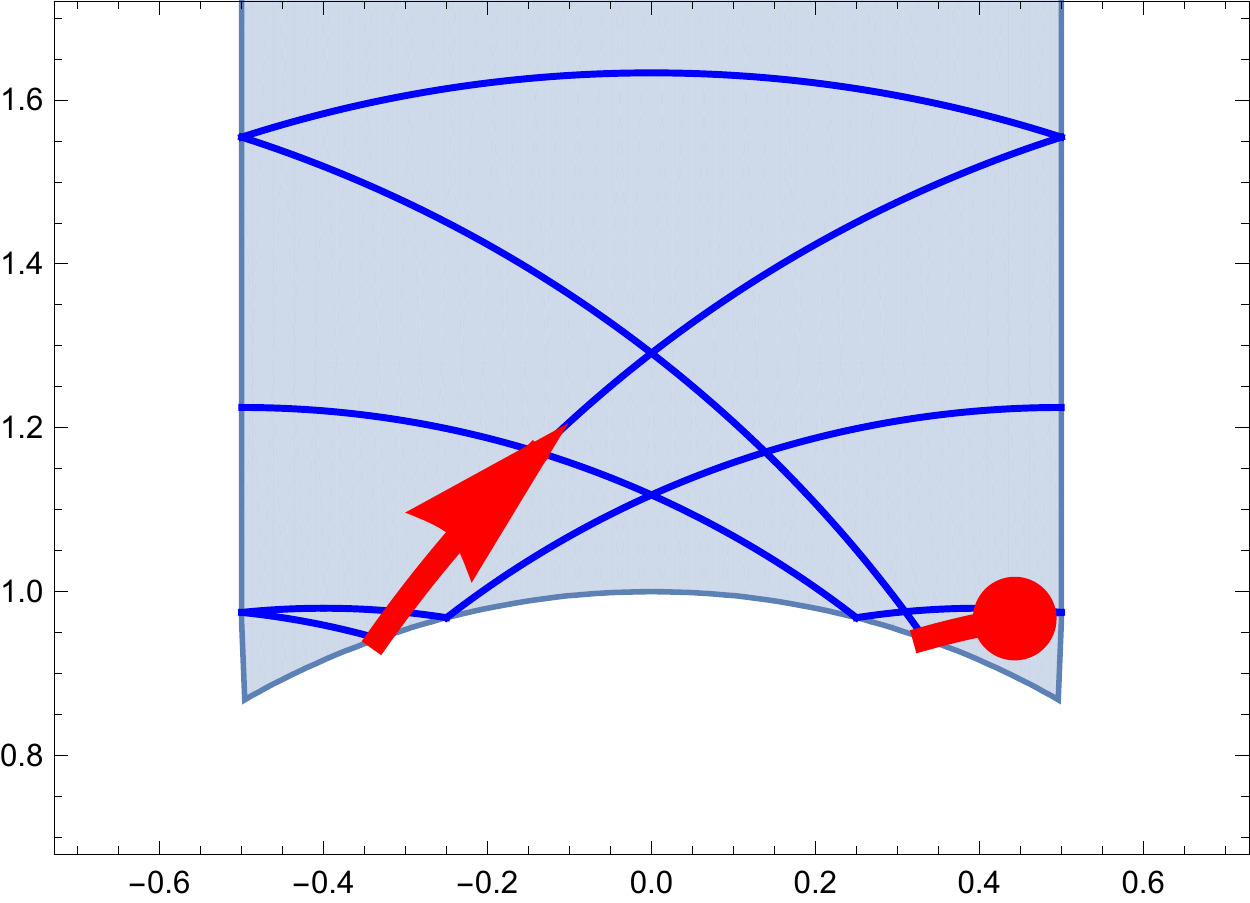}
                
{$L_{2}=S$}
        \end{subfigure}%
\quad        \begin{subfigure}[t]{0.3\textwidth}
                \centering
		\includegraphics[width=\textwidth]{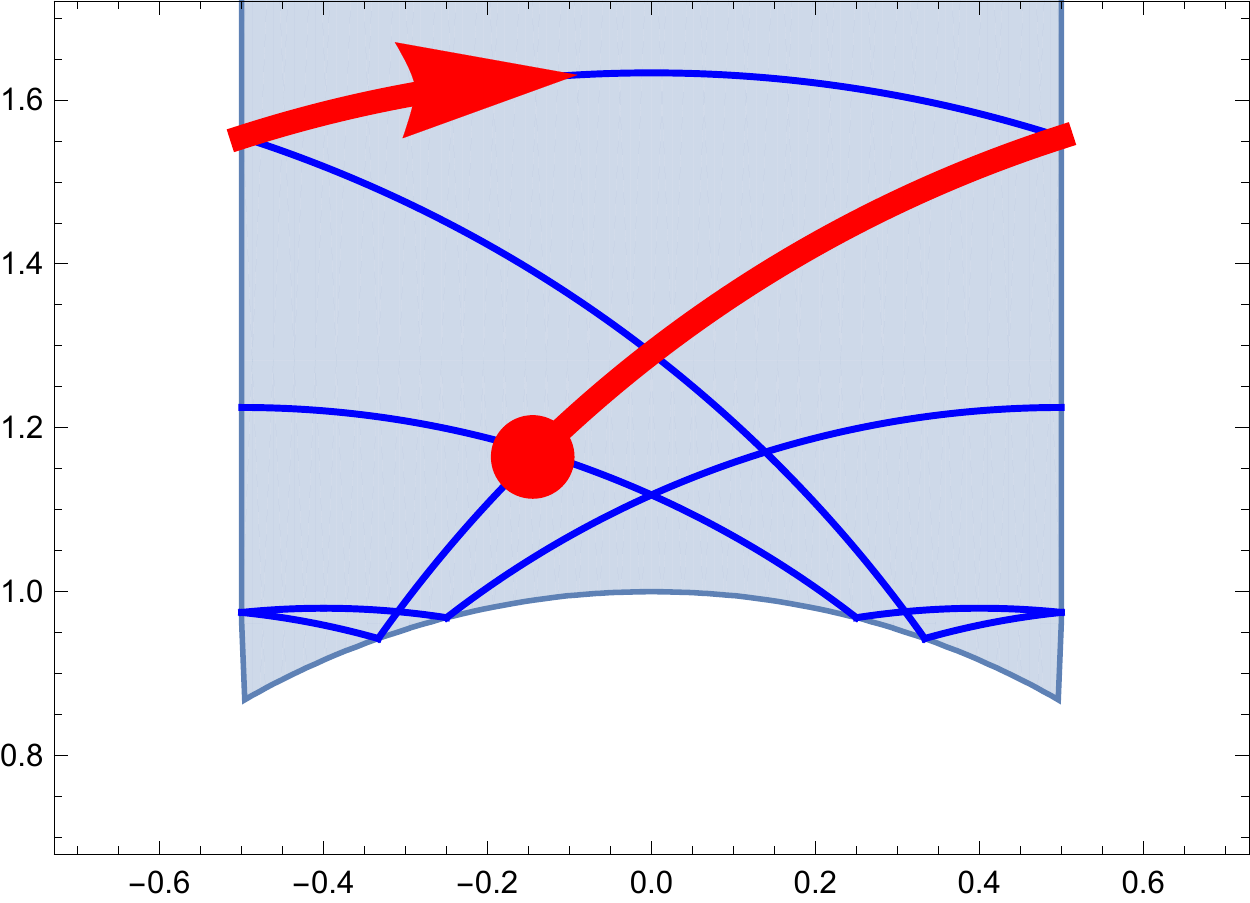}
                
{$L_{3}=T^{-1}$}
        \end{subfigure}%

\vspace{.1in}

        \begin{subfigure}[t]{0.3\textwidth}
                \centering
		\includegraphics[width=\textwidth]{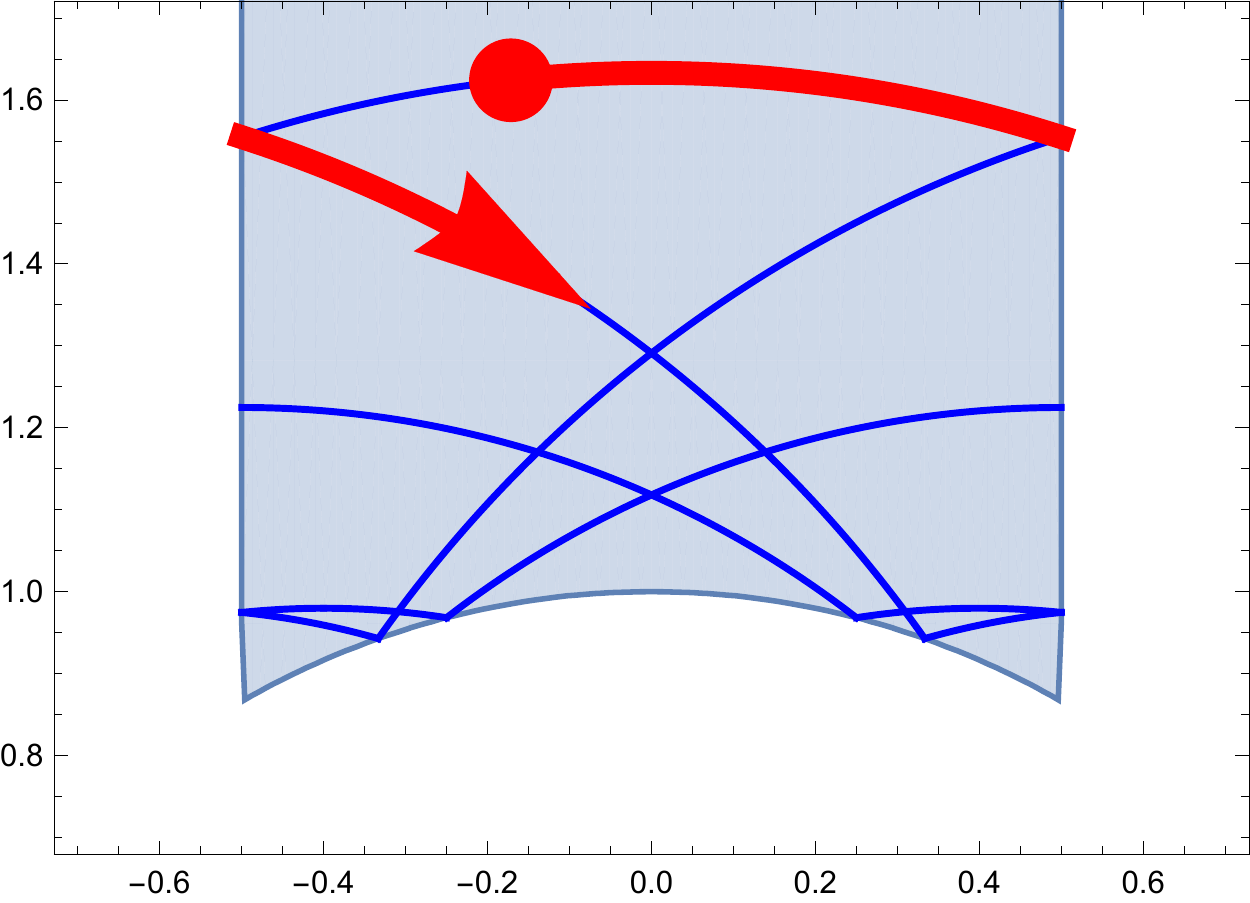}
                
{$L_{4}=T^{-1}$}
        \end{subfigure}%
\quad
        \begin{subfigure}[t]{0.3\textwidth}
                \centering
		\includegraphics[width=\textwidth]{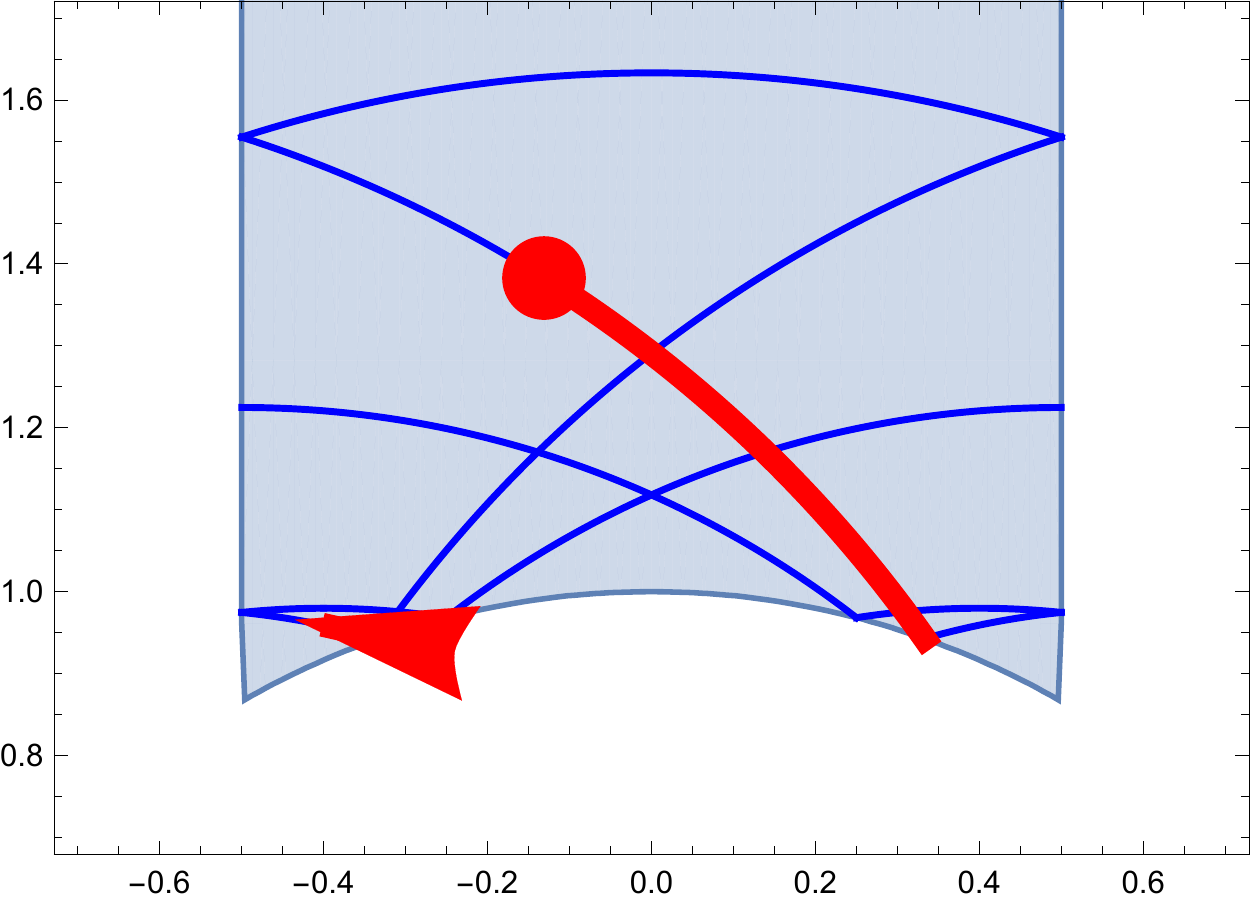}
                
{$L_{5}=S$}
        \end{subfigure}%
\quad        \begin{subfigure}[t]{0.3\textwidth}
                \centering
		\includegraphics[width=\textwidth]{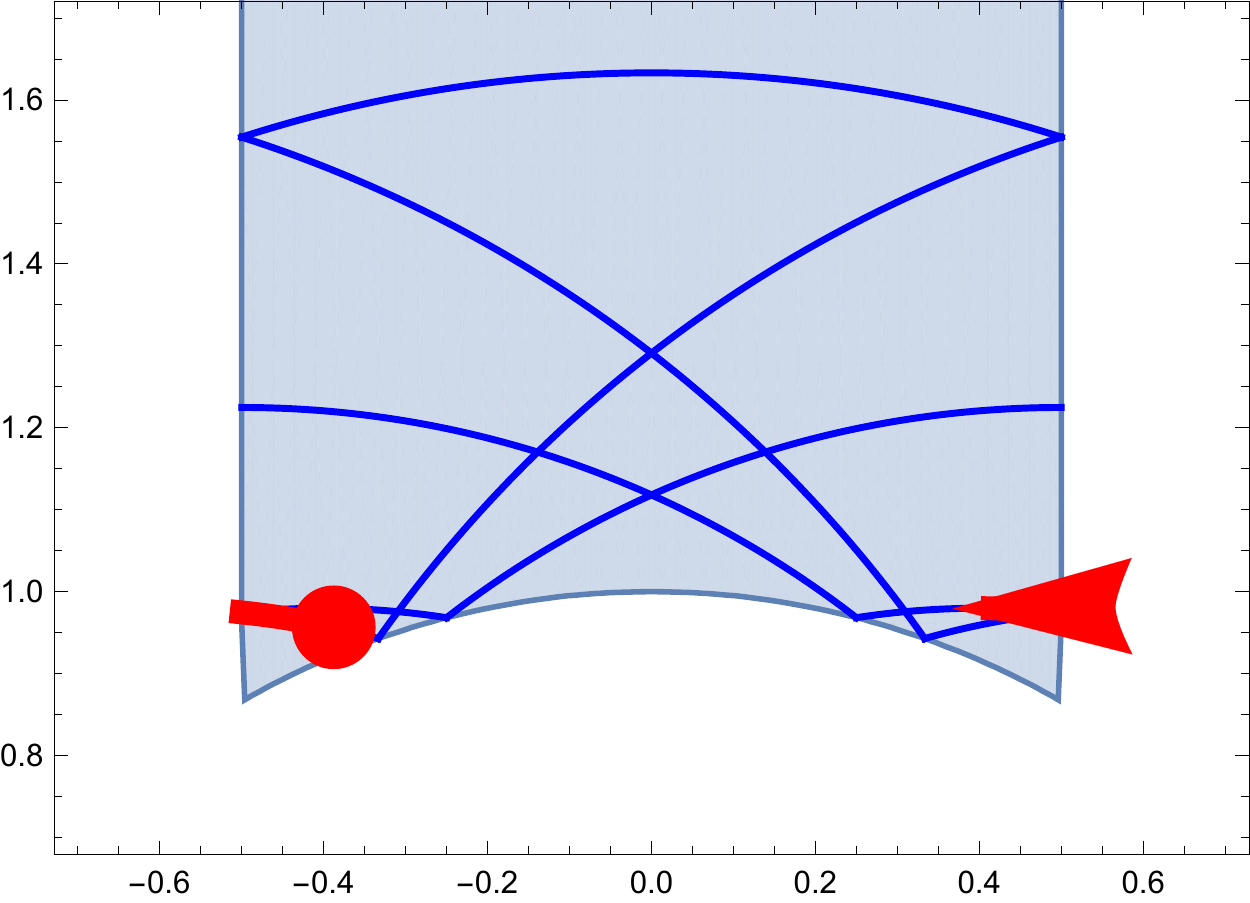}
                
{$L_{6}=T$}
        \end{subfigure}%

\vspace{.1in}

        \begin{subfigure}[t]{0.3\textwidth}
                \centering
		\includegraphics[width=\textwidth]{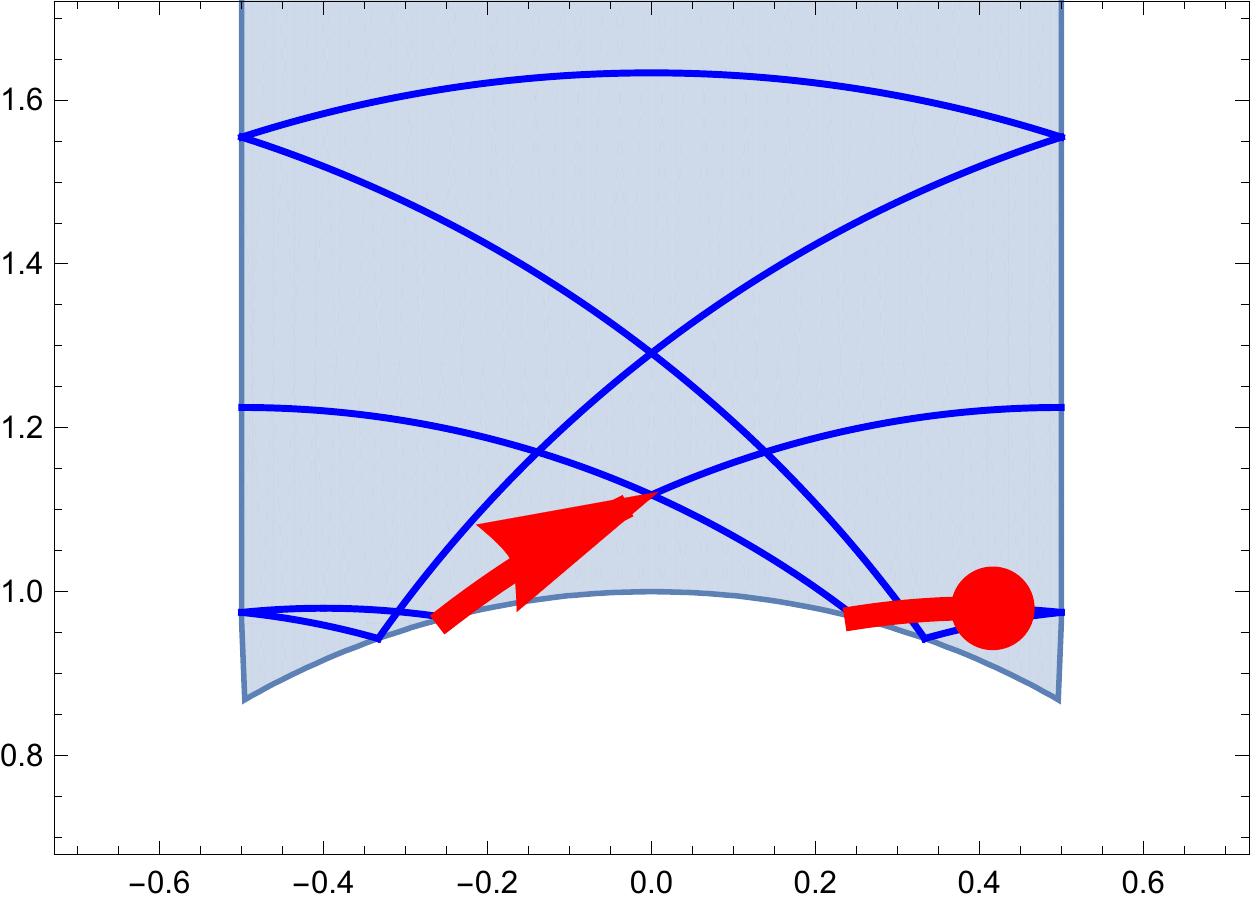}
                
{$L_{7}=S$}
        \end{subfigure}%
\quad
        \begin{subfigure}[t]{0.3\textwidth}
                \centering
		\includegraphics[width=\textwidth]{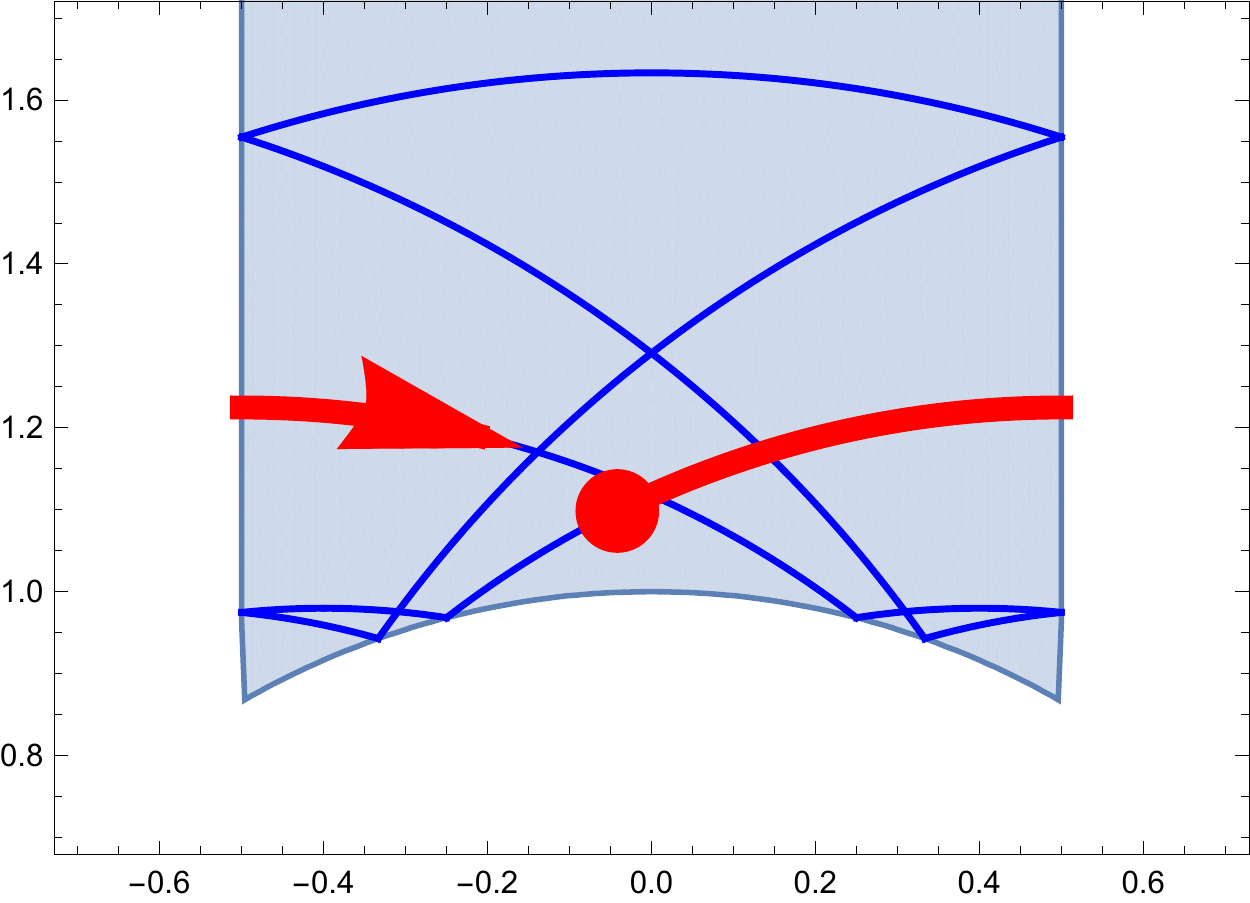}
                
{$L_{8}=T^{-1}$}
        \end{subfigure}%
\quad        \begin{subfigure}[t]{0.3\textwidth}
                \centering
		\includegraphics[width=\textwidth]{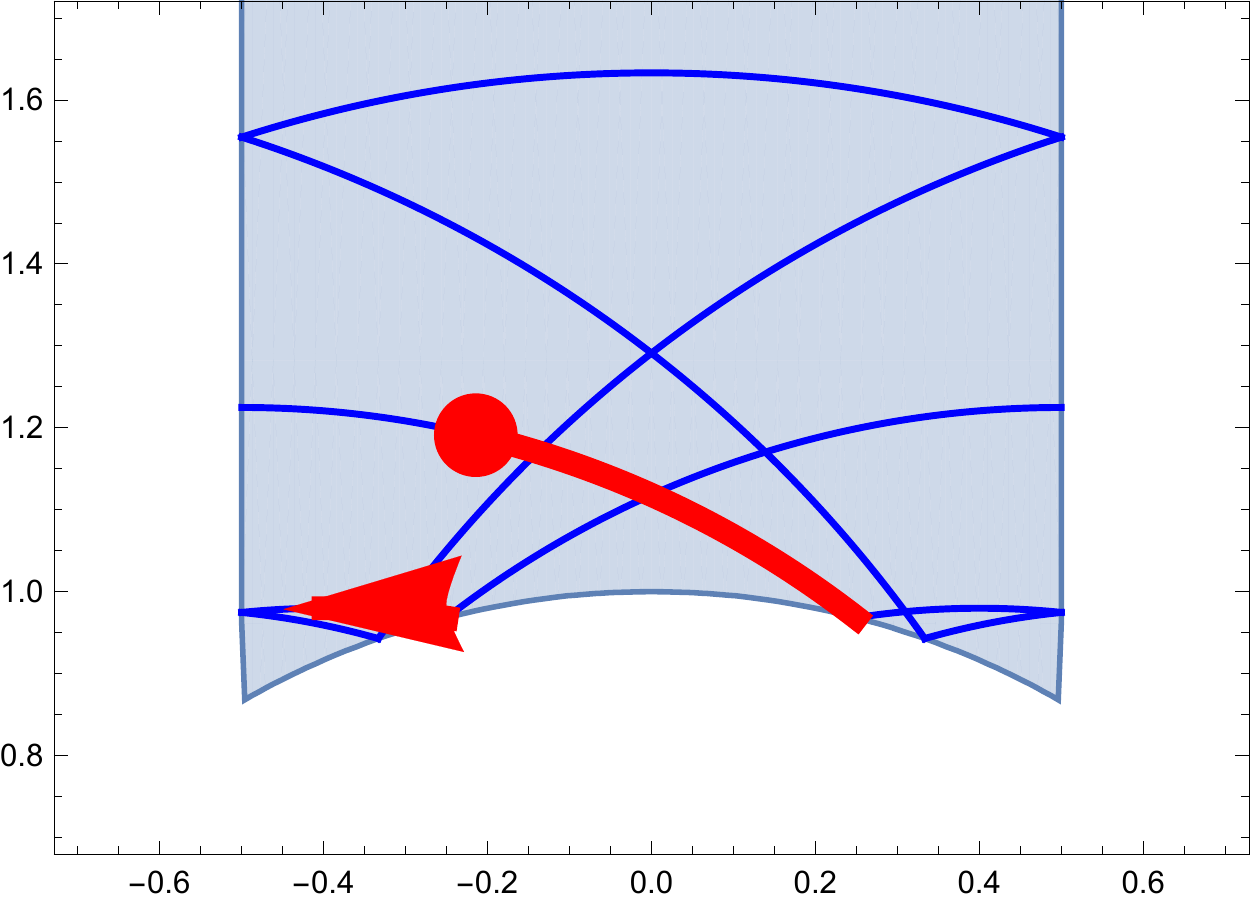}
                
{$L_{9}=S$}
        \end{subfigure}%

\caption{The cutting sequence of the geodesic flow} 
\label{AKfig:TSs}
\end{figure}

To illustrate this, consider again the example in \figref{AKfig:geodM}. The flow first hits the wall $\Re z=-1/2$, and must be translated by $L_{1}=T$ back inside $\cF$. Next the flow encounters the wall $|z|=1$, and is reflected using $L_{2}=S$. Continuing in this way (see \figref{AKfig:TSs}), we find that the cutting sequence of $(z,\gz)$ in 
\eqref{AKeq:zgzEg} is:
\be\label{AKeq:TSseq}
T,S,T^{-1},T^{-1},S,T,S,T^{-1},S, \dots,
\ee
repeating {\it ad infinitum}.
It is easy to see from the geometry that such sequences 
are 
some number of $T$'s or $T^{-1}$'s separated by single $S$'s. Computing these counts converts \eqref{AKeq:TSseq} into:
$$\underbrace{T}_{1},S,\underbrace{T^{-1},T^{-1}}_{2},S,\underbrace{T}_{1},S,\underbrace{T^{-1}}_{1},S, \dots,
$$
which corresponds to the sequence
\be\label{AKeq:TSto12s}
1,2,1,1,1,2,1,1,\dots
\ee
repeating.

It 
seems to have first been observed by Humbert \cite{Humbert1916}
that this sequence should be compared to the continued fraction expansion of the visual point $\ga$ of $(z,\gz)$. 
We write the 
continued fraction\index{continued fraction} expansion of any $x\in\R$ as
$$
x\ = \ a_{0}+\cfrac{1}{a_{1}+\cfrac{1}{a_{2} + \ddots}}=[a_{0},a_{1},a_{2},\dots],
$$
where $a_{0}\in\Z$ and the other $a_{j}\in\N$ are positive. These numbers are called the {\it partial quotients} of $x$, and we will sometimes call them ``digits'' or ``letters.'' 
For the visual point $\ga$ in \eqref{AKeq:gaIsEg}, we compute:
\be\label{AKeq:gaCFE}
\ga\ = \ {-7-2\sqrt6\over5} \ = \ [-3, 1, 1, \overline{1, 1, 1, 2}],
\ee
where the bar means repeating the last sequence of digits forever.

Comparing \eqref{AKeq:gaCFE} to \eqref{AKeq:TSto12s}, we see that the periodic parts match, up to cyclic permutation (since anyway a closed geodesic has no canonical ``starting'' point).
This leads us to the notion of a {\it reduced form} for $\ga$.

{\bf Definition:} 
A quadratic irrational $\ga$ is called {\it reduced} if it and its Galois conjugate $\overline \ga$ satisfy the
inequalities:
\be\label{AKeq:gaRed}
-1\ < \ \overline\ga \ < \ 0\ <\ 1\  <\  \ga.
\ee
A representative $M$ of a hyperbolic conjugacy class $[M]$ is also called {\it reduced} if its
visual point $\ga$
is.

{\bf Exercise:} A quadratic irrational $\ga$ is reduced iff its continued fraction is {\it exactly} (as opposed to eventually) periodic.

How should we reduce the representative $M$ in \eqref{AKeq:Mis}? It's actually quite easy. 
%
Note that, in general, if
 $\ga$ has continued fraction expansion
$$
\ga \ = \ [a_{0},\dots,a_{h},\overline{a_{h+1},\dots,a_{h+\ell}}],
$$
 then
\be\label{AKeq:cfeDigs}
\mattwo 011{-a_{0}}\cdot\ga\ = \ {1\over \ga-a_{0}}\ = \ [a_{1},\dots,a_{h},\overline{a_{h+1},\dots,a_{h+\ell}}].
\ee
That is,  such matrices eat away the first digit (this of course is the left-shift map from dynamics). 
For $\ga$ in \eqref{AKeq:gaCFE}, we 
could try
acting 
(on the left) by
$
\g_{0}=\mattwos 011{-1}\cdot\mattwos 011{-1}\cdot\mattwos 0113
$
to make an exactly periodic continued fraction $[\overline{1,1,1,2}]$. 
But
this matrix
$\g_{0}$
 has determinant $-1$, being an odd product of determinant $-1$ matrices. To act 
 instead
  by 
  an element of
  $\PSL_{2}$, we eat away one more digit, using the matrix
\be\label{AKeq:gamChVars}
\g \ = \ \mattwo 011{-1}\cdot\mattwo 011{-1}\cdot\mattwo 011{-1}\cdot\mattwo 0113 \ = \ 
\mattwo 25{-3}{-7},
\ee
to obtain
\be\label{AKeq:gaTil}
\widetilde \ga  \ = \ \g \ga \  = \ \frac{1+\sqrt{6}}{2} \ = \ [\overline{1,1,2,1}].
\ee
That is, we replace $M$ by
\be\label{AKeq:Mtil}
\widetilde M \  = \ \g M \g^{-1} \ = \
\mattwo 7543,
\ee
with $\widetilde M$ now reduced.

\

In this way, the geodesic flow corresponds to the symbolic dynamics\index{symbolic dynamics} on the continued fraction expansion of the visual point of the flow.
Actually, we have been quite sloppy; for example, it is not clear what 
one
should 
do if
the geodesic 
flow
passes through an elliptic point of the orbifold $\G\bk \bH$.
There is in fact a better way of encoding the cutting sequence, as elucidated 
beautifully
by Series \cite{Series1985}.
That said, the less precise (but more immediate) description given here will suffice for the purposes of our discussion below.

\

\

 \subsection{Binary Quadratic Forms}\label{sec:Gauss}\

 This will be a very 
 quick
 introduction 
 to
 an extremely well-studied and beautiful theory; see, e.g., \cite{Cassels1978} for a classical treatment.
The theory is largely due to  Gauss (building on Lagrange and Legendre), 
as developed
 his 1801
 magnum opus,
  {\it Disquisitiones Arithmeticae}.

For integers $A$, $B$, and $C$, let $Q=[A,B,C]$ denote the (integral) binary quadratic form\index{binary quadratic form} $Q(x,y)=Ax^{2}+Bxy+Cy^{2}$. 
The 
general problem being
addressed 
was: Given $Q$, what numbers does it represent? That is, for which numbers $n\in\Z$ do there exist $x,y\in\Z$ so that $Q(x,y)=n$?
The question is perhaps
inspired by the famous resolution 
in the 
case $Q=x^{2}+y^{2}$ by Fermat (and other special cases due to Euler and others).
See \cite{CoxBook} for a beautiful exposition of this problem.

Some observations: 

(i) If $A$, $B$, and $C$ have a factor in common, then so do all numbers represented by $Q$, and by dividing out this factor, we
may and will
assume henceforth that $Q$ is {\it primitive}, meaning $(A,B,C)=1$.

(ii) The set of numbers represented by $Q$ does not change if $Q$ is replaced by $Q'(x,y)=Q(ax+by,cx+dy)$, with $ad-bc=\pm1$;
this of course is nothing but an invertible
(over integers!)
 linear  change of variables. 

Gauss defined
 such a pair of
 forms $Q, Q'$ to be
 equivalent
 but
for   
  us it will be more convenient to use {\it strict} (some authors call this {\it proper}, or {\it narrow}) equivalence, meaning we only allow
  ``orientation-preserving'' transformations. That is, we will write $Q\sim Q'$ only when there is some 
   $\g\in \SL_{2}(\Z)$ (as opposed to $\GL_{2}$) with $Q=Q'\circ\g$.

{\bf Exercise:}
This is indeed an {\it action}, that is, $ (Q\circ \g_{1})\circ\g_{2}=Q\circ(\g_{1}\g_{2})$, and hence $\sim$ is an equivalence relation. 

For a given $Q$, the set of all $Q'\sim Q$ is called a {\it class} (or equivalence class) and denoted $[Q]$.
Because we are 
considering
strict equivalence here, 
this is 
often called 
the  {\it narrow} class of $Q$.

{\bf Exercise:}
If $Q\sim Q'$ then $D_{Q}=D_{Q'}$, where 
\be\label{AKeq:discr}
D_{Q}:=\discr(Q)=B^{2}-4AC
\ee 
is the discriminant. That is, the discriminant is a class function (invariant under equivalence). 
Observe 
 that discriminants
  are quadratic residues $\pmod 4$, and hence 
  $D_{Q}\equiv0$ or $1\pmod 4$.

{\bf Exercise:}
When $D<0$, the form $Q$ is {\it definite}, that is, it only takes either positive or negative values, but not both. 
When $D>0$, the form is {\it indefinite}, representing both positive and negative numbers.

If $\discr(Q)=
D=0,$ 
%
or more generally, if
 $D=D_{0}^{2}$ is a perfect square, then $Q$ is the product of two linear forms. 
 Then the representation question is much less interesting, 
 and will be left as an exercise.
We
 exclude this case
going forward.

{\bf Exercise:} Let
\be\label{AKeq:gaQ}
\ga_{Q} \ = \ {-B+\sqrt{D}\over 2A}
\ee
be the root  of $Q(x,1)$ (assuming $A\neq0$), and suppose $Q'= Q\circ \g$. Then $\ga_{Q'}=\g^{-1}\cdot\ga_{Q}$, where 
the action here of
$\g^{-1}$ is by fractional linear transformations.

We have seen that if two forms are equivalent, then their discriminants agree.
It is then natural to ponder about the converse:
Does $\discr(Q)=\discr(Q')$ imply that $Q\sim Q'$?

To study this question, let $\sC_{D}$ be the set of all inequivalent, primitive classes 
having discriminant $D$, 
$$
\sC_{D}:=\{[Q]:\discr(Q)=D\},
$$
and 
let $h_{D}:=|\sC_{D}|$ be 
its size;
this is called the (narrow) {\it class number}. 

This set $\sC_{D}$ is now called the {\it class group} (it turns out there is a composition process under which $\sC_{D}$ inherits the  structure of an abelian group, but this fact will not be needed for our investigations; for us, $\sC_{D}$ is just a set). 
If having the same discriminant implied equivalence, then all class numbers would be $1$.
This turns out to be false,
but actually it is not off by ``very'' much, in the
following sense:
the class number\index{class number}
is always finite.

\begin{theorem}[Gauss]\label{AKthm:hD}
For any non-square integer $D\equiv0,1(4)$, we have: 
$$
1\ \le \  h(D) \ <  \ \infty.
$$
\end{theorem}

It is 
easy to see that
 $h(D)\ge1$.
  Indeed, if $D\equiv0(4)$, then $x^{2}-\frac D4 y^{2}$ is primitive of discriminant $D$; and if $D\equiv1(4)$, then $x^{2}+xy-{D-1\over 4}y^{2}$ works.

We give a very quick sketch of the very well-known \thmref{AKthm:hD}, as it will be relevant to what follows. The proof decomposes according to whether or not the class is definite.

\pf[Sketch when $D<0$.]
Let $\ga_{Q}$ be the root of $Q(x,1)$ as in \eqref{AKeq:gaQ}; since $Q$ is definite, $\ga_{Q}\in\bH$.
We already discussed the standard fundamental domain $\cF$ for $\G\bk \bH$, so we know
 (from the exercise below \eqref{AKeq:gaQ}) 
 that there is a transformation
$\g\in\SL_{2}(\Z)$ taking $\ga_{Q}$ to $\cF$. 
Such a transformation $\g$ is moreover unique (up to technicalities when $\ga_{Q}$ is on the boundary of $\cF$),
and gives rise to a unique {\it reduced} representative $Q'= Q\circ\g^{-1}$ in the class $[Q]$.
It remains to show there are only finitely many such having a given discriminant $D<0$. 
This
is an easy exercise using $|\Re(\ga_{Q})|\le 1/2$ and $|\ga_{Q}|\ge1$;
indeed, one finds that
$$
A\ \le\ \sqrt{|D|\over 3},
\quad\text{and}\quad |B| \ \le \  A,
$$
and hence the number of reduced forms is finite.
%
\epf

{\bf Example:} Take $D=-5$. This $D$ is not congruent to $0$ or $1$ mod $4$, so instead we consider $D=-20$. Then $A\le\sqrt{20/3}$, that is, $A$ is at most $2$. If $A=1$, then $B=0$ or $\pm1$. The latter case,
$B=\pm1$, gives
 no integral solutions to $C=(B^{2}-D)/(4A),$ but the former gives $C=5$, corresponding to the reduced form $Q_{0}=x^{2}+5y^{2}$.
Next we consider $A=2$, whence $B\in[-2,2]$. Only the choices $B=\pm2$ give integral values for $C=3$, corresponding to forms $Q_{\pm}=2x^{2}\pm2xy+3y^{2}$. Actually this turns out to be a boundary case, and the two forms $Q_{\pm}$ are equivalent. Hence the class group $\sC_{-20}=\{Q_{0},Q_{+}\}$ has class number  $h_{-20}=2$. 
The fact that 
this
class number is not $1$ is well known to be related to the failure of unique factorization in $\Z[\sqrt{-5}]$, e.g., the number $6$  factors both as $6=2\cdot3$ and as $6=(1+\sqrt5i)(1-\sqrt5i)$.

\pf[Sketch when $D>0$.]
Again we assume $A>0$.
This case is more subtle, 
as
the root $\ga_{Q}$ and its Galois conjugate are real, so 
cannot
 be moved to 
  the fundamental domain $\cF$.
   Instead we notice that, since $\ga_{Q}$ is real and quadratic, it has an eventually periodic continued fraction expansion, and transformations $\g\cdot\ga_{Q}$ simply add (or subtract) letters to (or from) this expansion. In particular, there is a $\g\in\SL_{2}(\Z)$ which makes the continued fraction exactly periodic; that is, 
the transformed root
$\ga_{Q\cdot\g}$ 
then satisfies
the
familiar
 condition \eqref{AKeq:gaRed}. (This is our first hint that indefinite forms are related to closed geodesics!) We thus call an indefinite form $Q$ {\it reduced} if its root $\ga_{Q}$ is reduced. 
 As before, it is easy to see using
 $$
\ga_{Q} \ = \ 
{-B+\sqrt{D}\over 2A} \ > \ 1,\qquad -1\ < \ \overline\ga_{Q} \ = \ {-B-\sqrt{D}\over 2A} \ < \ 0,
 $$
 that
 $$
0\ < \ -B\ < \ \sqrt{D},\qquad 
\frac12(\sqrt D+B) \ < \ A \ < \ \frac12(\sqrt D -B),
 $$
 which forces the class number to be finite. 
\epf

{\bf Example:} Take $D=7,$ or rather $D=28$, since $7\not\equiv0,1(4)$. Then the possibilities for $B$   range from $-1$ to $-5$. To solve $B^{2}-D=4AC,$ we must have $B^{2}-D\equiv0(4)$, which leaves only $B=-2$ or $B=-4$. In the former case, the possible positive divisors of $(B^{2}-D)/4=-6$ are $A=1$, $2$, $3$, or $6$, of which only $A=2$ and $3$ lie in the range $\frac12(\sqrt{28}-2)<A<\frac12(\sqrt{28}+2)$.
These two give rise to the forms 
\be\label{AKeq:Q1Q2}
Q_{1}=2x^{2}-2xy-3y^{2}
\qquad\text{and }\qquad
Q_{2}=3x^{2}-2xy-2y^{2}.
\ee 
The latter case of $B=-4$ leads in a similar way to the two reduced forms $Q_{3}=x^{2}-4xy-3y^{2}$ and $Q_{4}=3x^{2}-4xy-y^{2}$. 

Note that, unlike the definite case, reduced forms are {\it not} unique in their class, as any even-length cyclic permutation of the continued fraction expansion of $\ga_{Q}$ gives rise to another reduced form. So we cannot conclude from the above computation that $h_{28}\overset{?}{=}4$. Writing $\ga_{j}$ for the larger root of $Q_{j}$, $j=1,\dots,4$, we compute the continued fractions:
\be\label{AKeq:gaCFEs}
\ga_{1}=[\overline{1,1,4,1}],\ 
\ga_{2}=[\overline{1,4,1,1}],\ 
\ga_{3}=[\overline{4,1,1,1}],\ 
\ga_{4}=[\overline{1,1,1,4}].
\ee
It is then easy to see by inspection that $Q_{1}$ and $Q_{3}$ are equivalent, as are $Q_{2}$ and $Q_{4}$, e.g., 
$$
Q_{1}\circ\left[\mattwos011{-1}\mattwos011{-1}\right]^{-1}=Q_{3}.
$$ 
({\bf Exercise:} Verify this and compute the change of basis matrix to go from $Q_{2}$ to $Q_{4}$.) Thus $h_{28}=2$. 
\\

Class groups and class numbers are extremely mysterious.
A discriminant is defined to be {\it fundamental} if it is the discriminant of a (
quadratic, in our context) 
field; such $D$'s are either $\equiv1\pmod 4$ and squarefree, or  divisible by $4$ with $D/4$ squarefree and $\equiv2,3\pmod4$.
Dirichlet's Class Number Formula\index{Dirichlet's Class Number Formula} (see, e.g., \cite{DavenportBook, IwaniecKowalski}) gives 
an approach 
to studying class groups: If $D$ is a fundamental discriminant\index{fundamental discriminant}, then
\be\label{AKeq:DCNF}
h_{D} \ = \  \sqrt {|D|}\,L(1,\chi_{D})  \times\twocase{}{1/(2\pi),}{if $D\le-5$, or}{1/\log\gep_{D},}{if $D>0$.}
\ee
Here $L(1,\chi_{D})$ is called a ``special $L$-value,'' 
and for $D>0$, the 
factor
 $\gep_{D}\in\Q(\sqrt D)$ is 
 determined by: 
\be\label{AKeq:gepDdef}
\gep_{D} \ = \ {t+s\sqrt D\over 2},
\ee
where $(T,S)=(t,s)$ is the least solution to the Pellian equation\index{Pell equation} 
\be\label{AKeq:Pell}
T^{2}-S^{2}D=4.
\ee
It follows from \eqref{AKeq:Pell} that $\gep_{D}$ is a unit in $\Q(\sqrt D)$, as its algebraic norm is $
\gep_{D}\overline\gep_{D}=1$. (If the Pell equation \eqref{AKeq:Pell} has no solutions with $4$ replaced by $-4$ on the right side, then $\gep_{D}$ is the fundamental unit;\index{fundamental unit} otherwise it is the square of the latter.)
The $L$-value, which we will not bother to define (as it is not relevant to our discussion),
 is so fascinating an object about which one could say so much,
  that
   we will
instead
 say 
very little.
For example, it is not hard to show that 
 $$
 L(1,\chi_{D}) \ \le 
 C \log |D|.
 $$
 Siegel famously proved  \cite{Siegel1935, Landau1935} 
 the reverse inequality, that
 for any $\vep>0$,  
$$
L(1,\chi_{D})   \ \ge \ C_{\gep} \cdot 
|D|^{-\vep}
,
$$ 
though the constant $C_{\vep}$ is {\it ineffective}, that is, the proof does not give any means to estimate 
it
for any given $\gep$. 
(These days, we have other, not as strong, but 
on the other hand,
effective, estimates, see \cite{Goldfeld1976, GrossZagier1986}.)
Either way, we may think of the $L$-value as very roughly being of size $1$.
Then definite class groups are, very roughly, of size
$$
h_{-D} \ \approx \ \sqrt{|D|},
$$
while indefinite ones are of size
\be\label{AKeq:hDis}
h_{D} \ \approx \ \sqrt{D}/\log \gep_{D}.
\ee
(We will not give the symbol $\approx$ a precise meaning here.)
As a consequence, 
we obtain Gauss's conjecture, that for definite class numbers, 
$
h_{-D} 
 \to  \infty
 $
,
as
$
-D \to  -\infty
,
$
see \cite{Deuring1933, Heilbronn1934}.

For indefinite forms, the behavior is even more mysterious, and  it is a longstanding conjecture that infinitely 
often the 
class number is $1$:
\be\label{AKeq:CNO}
\liminf
h_{D} \ \overset{?}{=} \ 1
,
\ee
where the limit is over fundamental $D\to+\infty$. (If one does not require fundamental discriminants, this problem was apparently solved 
long ago
 by Dirichlet, see \cite{Lagarias1980}.)
In light of \eqref{AKeq:hDis}, this conjecture suggests that the 
unit 
$\gep_{D}$, defined in \eqref{AKeq:gepDdef}, should infinitely often be massive, of size about $e^{\sqrt D}$. Today no  methods are known to force such large solutions to the Pell equation \eqref{AKeq:Pell}, despite
rather convincing
 evidence (see, e.g., \cite{CohenLenstra1983}) that this event is far from rare. On the other hand, it is quite easy to make giant class numbers, since one can force very small units,\index{fundamental unit} e.g., by taking $D$'s of the form $D=t^{2}-4$, one sees that $\gep_{D}={t+\sqrt D\over 2}\approx \sqrt D$, and $h_{D}\approx \sqrt D/\log D$ is as large as possible. Another long-standing conjecture in the indefinite case is that the average class number is roughly bounded, in the crude (more refined conjectures are available) sense that:
\be\label{AKeq:avgHD}
\sum_{0<D<X}h_{D} \ \overset{?}{=} \ X^{1+o(1)}.
\ee
We have insufficient space to delve further into this fascinating story, so leave it there.

\subsection{Duke's Theorem}\label{sec:Duke}

We now combine the previous two sections, giving a bijection between primitive, indefinite classes $[Q]$ and primitive, oriented, closed geodesics $\g$ on the modular surface. 
This equivalence was apparently first observed by Fricke and Klein \cite{FrickeKlein1890}, and is also discussed in many places, e.g., \cite{Cassels1978, SarnakThesis, Hejhal1983, IwaniecKowalski, Sarnak2007}.

We first attach a form $Q$ to a given hyperbolic matrix
$M=\mattwos abcd$, 
the  passage being through equating $\ga$'s in \eqref{AKeq:gaIs} and \eqref{AKeq:gaQ},
as follows.
Pattern 
matching \eqref{AKeq:gaIs} with \eqref{AKeq:gaQ}, we obtain a preliminary (possibly imprimitive) set of variables $B_{0}=d-a$, $A_{0}=c$, $D_{0}=\tr^{2}M-4$, leading to 
$$
C_{0} \ = \ {B_{0}^{2}-D_{0}\over 4A_{0}} \ = \ -b.
$$
Setting $s=\gcd(A_{0},B_{0},C_{0})$, the primitive quadratic form $[A_{0},B_{0},C_{0}]/s$ is almost what we want, but does not quite work, because $M$ is really in $\PSL_{2}(\Z)$, and as is, $-M$ could give a
 different form.
 To fix this, we set
\be\label{AKeq:MtoQ}
Q \ = \ {\sgn(\tr M)\over s}[c,d-a,-b],
\ee
which is now well defined on $\PSL_{2}(\Z)$. The discriminant of this form is
\be\label{AKeq:DiscrIs}
D\ = \ {\tr^{2}M-4\over \gcd(c,d-a,b)^{2}}.
\ee
For example, if we take $\widetilde M=\mattwos7543$ in \eqref{AKeq:Mtil},
then $s=\gcd(4,-4,-5)=1$ and $\sgn(\tr M)=+1$, so $\widetilde M$ corresponds to
 the binary quadratic form  $\widetilde Q=4x^{2}-4xy-5y^{2}$ of discriminant $96$. 
Note that if we had done the same operation starting with $M$ in \eqref{AKeq:Mis},
 the 
 form
 $Q=-5x^{2}-14xy-5y^{2}$, also of discriminant $96$, and of course, $Q=\widetilde Q\circ \g^{-1}$ with change of variables matrix $\g$  given by \eqref{AKeq:gamChVars}. ({\bf Exercise:} Verify all this.)

To invert the map \eqref{AKeq:MtoQ} given some $Q=[A,B,C]$ with discriminant $D>0$ and not a perfect square, we seek a matrix $M=\mattwos a{-Cs}{As}d\in\SL_{2}(\Z)$ so that  $a+d>0$, say, and $d-a=Bs$. Inserting the last identity into the determinant equation and completing the square gives:
$$
1 \ = \ 
a^{2}+Bsa+ACs^{2} 
\ = \ 
\frac14\left(2a+{Bs}\right)^{2}-{1\over 4}Ds^{2} 
.
$$
Multiplying both sides by $4$, we come to the familiar Pellian equation \eqref{AKeq:Pell}; if $(t,s)$ is a fundamental solution, then
$
a =  {(t-Bs)
/ 2},
$ 
$d  =  {(t+Bs)
/2
},
$
 and we have found our desired hyperbolic matrix
$$
M \  = \ 
\mattwo{{(t-Bs)
/ 2}}{-Cs}{As}{{(t+Bs)
/ 2}}.
$$
That $M$ is primitive follows from the fundamentality of the solution $(t,s)$ ({\bf Exercise}). 
For example, 
to turn $Q_{1}=2x^{2}-2xy-3y^{2}$ of discriminant $D=28$ into a closed geodesic, we find the fundamental solution $(t,s)=(16,3)$ to \eqref{AKeq:Pell}, leading to
$
M_{1}=\mattwos{11}965.
$

Note that the key to finding $M$ above is to solve a Pellian equation, which itself goes through  continued fractions; here is a more direct
version of this inverse map. 
We first state the following simple

{\bf Exercise:} The matrix 
\be\label{AKeq:MfromAs}
M \ = \ \mattwo{a_{0}}110\cdot\mattwo{a_{1}}110\cdots\mattwo{a_{\ell}}110
\ee
fixes the real quadratic irrational $\ga$ having continued fraction expansion $\ga=[\overline{a_{0},a_{1},\dots,a_{\ell}}]$. [Hint: Compare to \eqref{AKeq:cfeDigs}.]

Recall from 
 \eqref{AKeq:MgaEqga} 
 that the desired matrix $M_{1}$
  fixes $\ga_{1}$, where $\ga_{1}=(2+\sqrt{28})/4$ is the root of $Q_{1}(x,1)$.
As the continued fraction expansion of $\ga_{1}$ is given  in \eqref{AKeq:gaCFEs}, it is trivial to find the corresponding matrix $M_{1}$:
$$
M_{1} \ = \ \mattwo{1}110\cdot\mattwo{1}110\cdot\mattwo{4}110\cdot\mattwo{1}110 \ = \ \mattwo{11}965,
$$
again.
If the product in \eqref{AKeq:MfromAs} had odd length, one would obtain a matrix of determinant $-1$, whose square is the desired 
element of $\SL_{2}(\Z)$. This corresponds to the situation when  $T^{2}-Ds^{2}=-4$ is solvable;
then, writing $(t,s)$ for the least solution,
we find that $\gep_{D}$ is expressed as
 $\eta_{D}^{2}=\gep_{D}$,
where $\eta_{D}=(t+s\sqrt D)/2$ is the fundamental unit.\index{fundamental unit}
\\

{\bf Definition:} The {\it discriminant} of a closed geodesic\index{discriminant of a closed geodesic} $\g$ on the modular surface, or its corresponding hyperbolic conjugacy class, is defined to be that of its associated equivalence class of binary quadratic forms.
Explicitly, the discriminant of $M=\mattwos abcd$ is $D=(\tr^{2}M-4)/s^{2}$, where $s=\gcd(c,d-a,b)$.
\\

The idea in Duke's Theorem
 is to look at the equidistribution of closed geodesics on the modular surface, but instead of studying them individually, one may, like classes of binary forms, group them by discriminant. First some examples.

 \begin{figure}
        \begin{subfigure}[t]{0.3\textwidth}
                \centering
		\includegraphics[width=\textwidth]{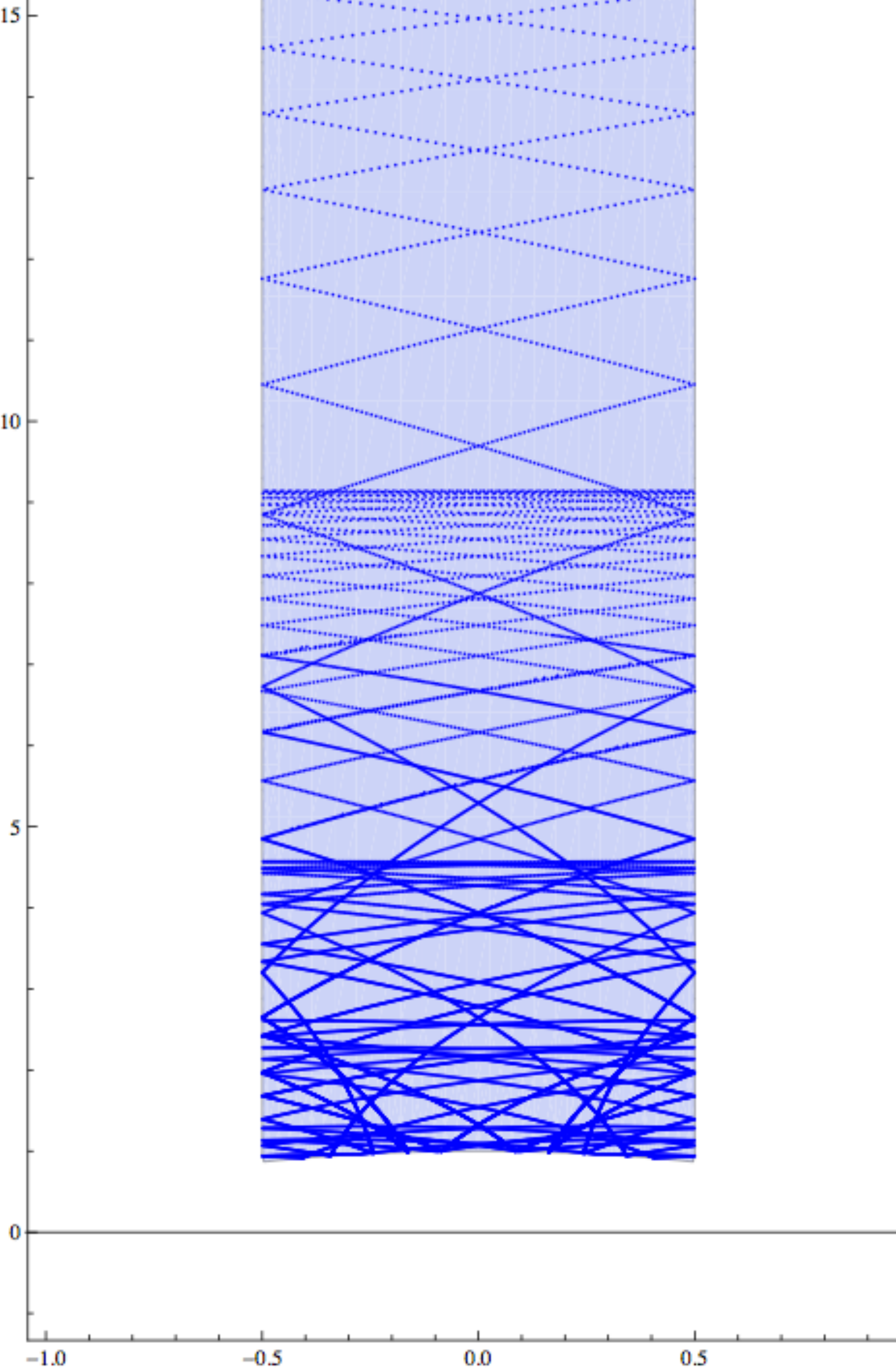}
                \caption{$D=1337$}
                \label{AKfig:D1337}
        \end{subfigure}%
\quad
        \begin{subfigure}[t]{0.3\textwidth}
                \centering
		\includegraphics[width=\textwidth]{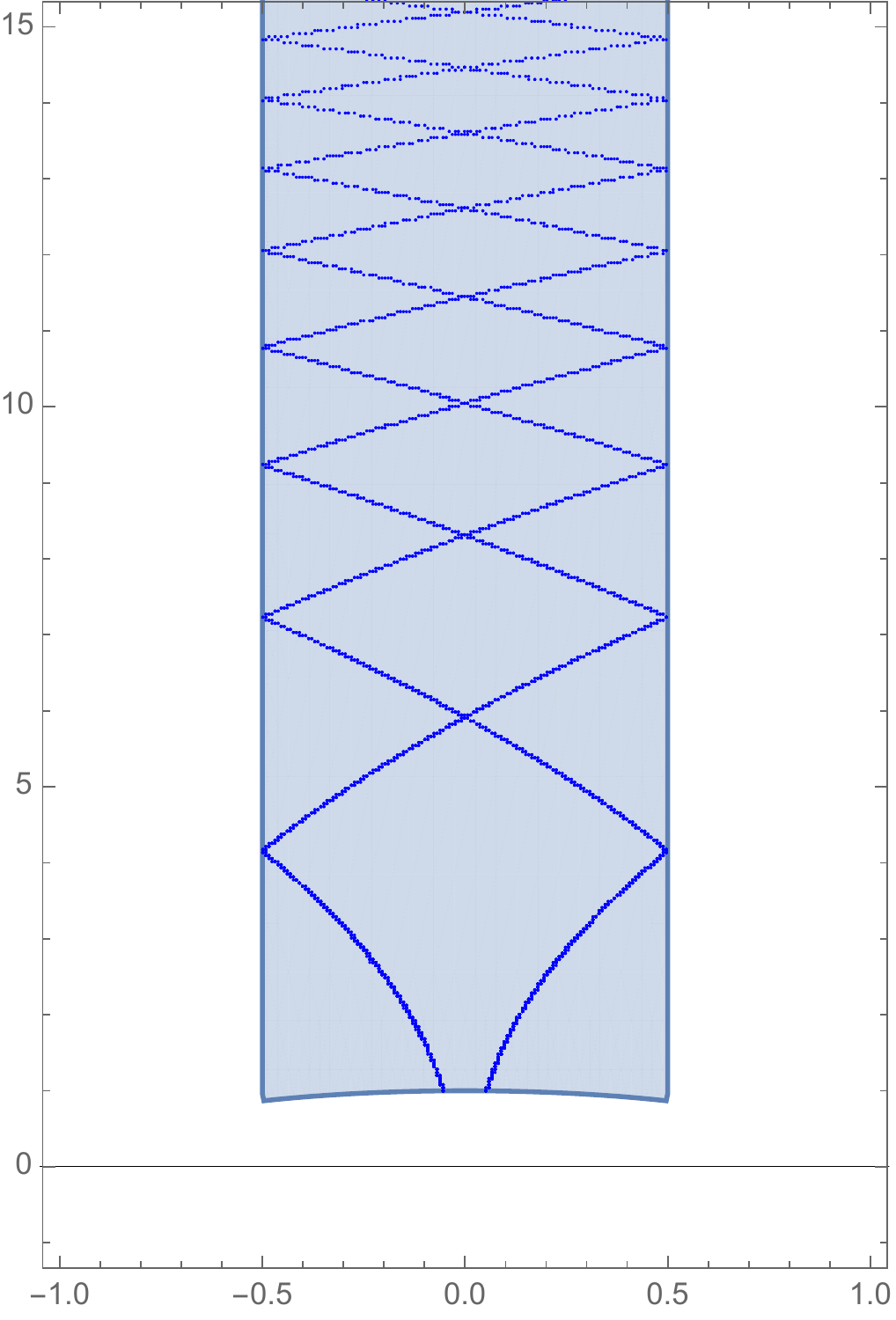}
                \caption{$D=1365:\ \ga_{1},\ga_{2}$} 
                \label{AKfig:D136512}
        \end{subfigure}
\quad
        \begin{subfigure}[t]{0.3\textwidth}
                \centering
		\includegraphics[width=\textwidth]{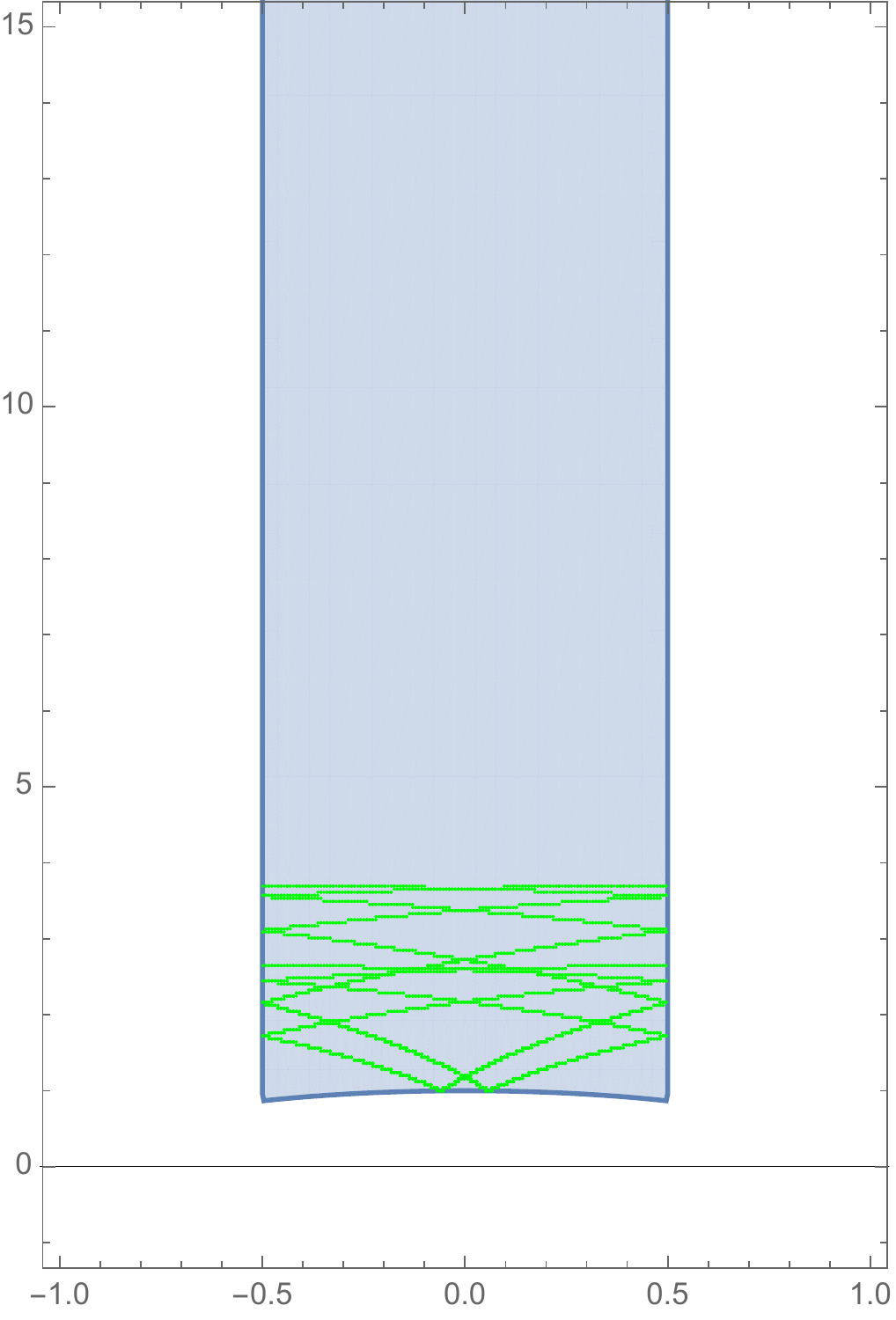}
                \caption{$D=1365:\ \ga_{3},\ga_{4}$} 
                \label{AKfig:D136534}
        \end{subfigure}
\quad
        \begin{subfigure}[t]{0.3\textwidth}
                \centering
		\includegraphics[width=\textwidth]{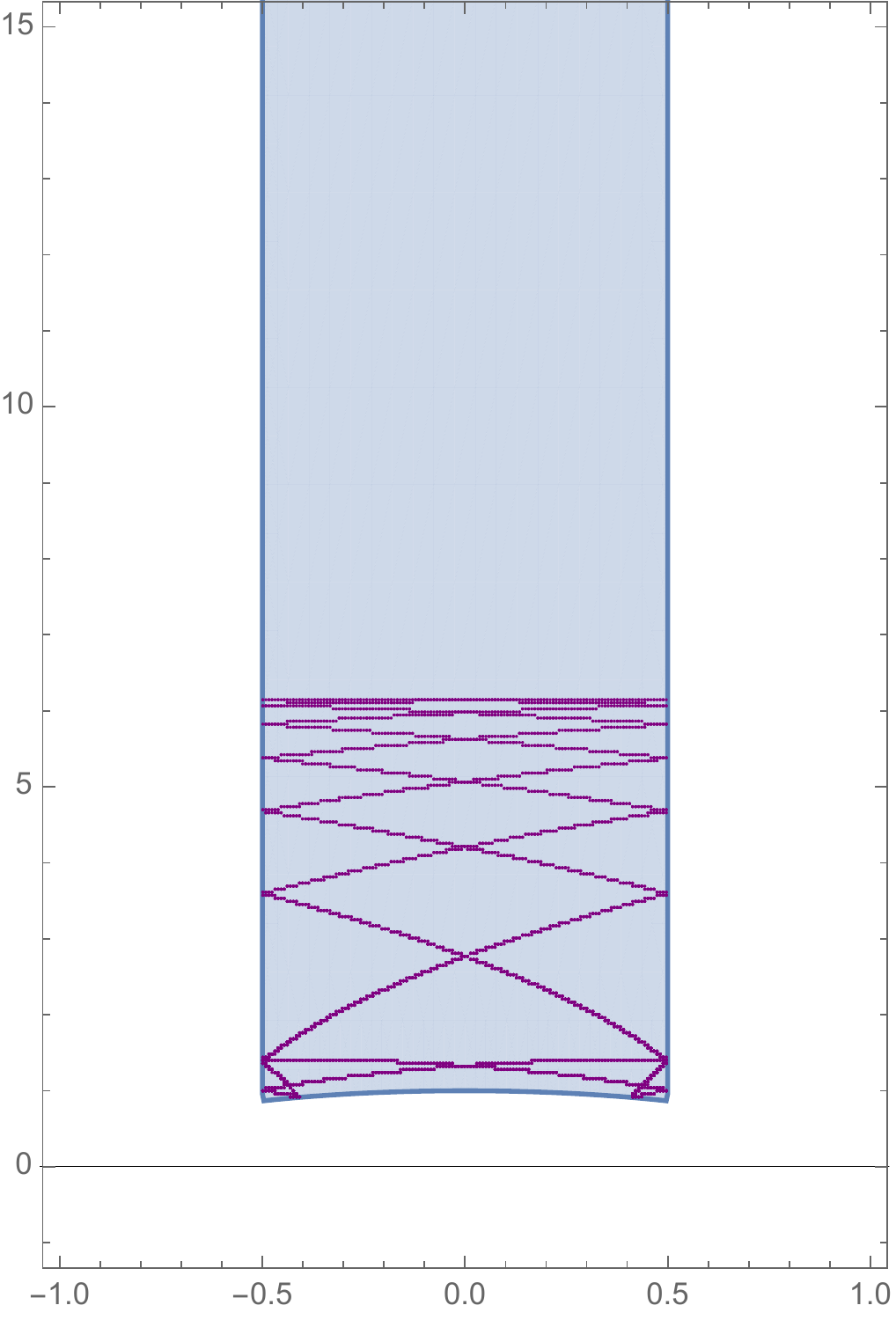}
                \caption{$D=1365:\ \ga_{5},\ga_{6}$} 
                \label{AKfig:D136556}
        \end{subfigure}
\quad
        \begin{subfigure}[t]{0.3\textwidth}
                \centering
		\includegraphics[width=\textwidth]{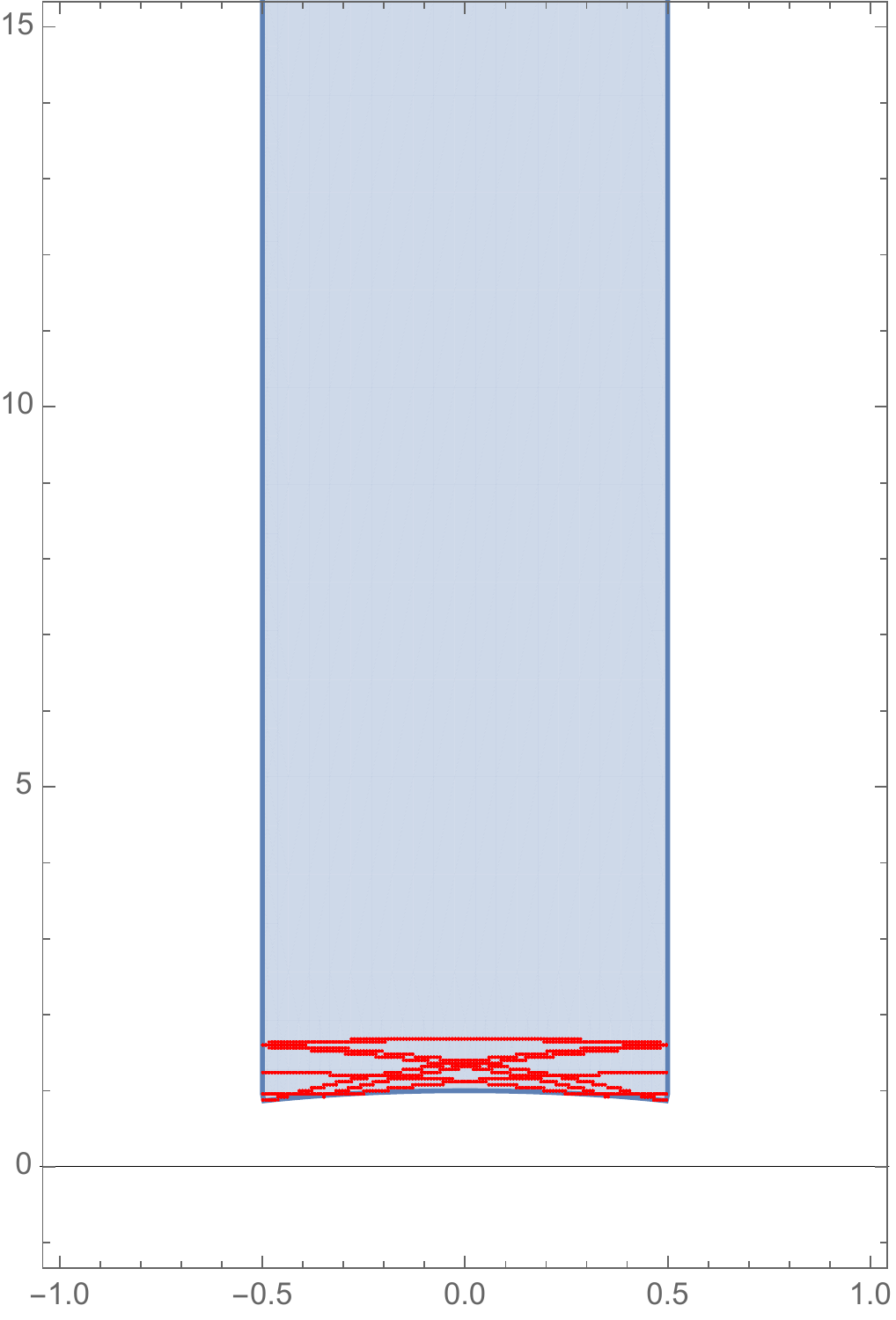}
                \caption{$D=1365:\ \ga_{7},\ga_{8}$} 
                \label{AKfig:D136578}
        \end{subfigure}
\quad
        \begin{subfigure}[t]{0.3\textwidth}
                \centering
		\includegraphics[width=\textwidth]{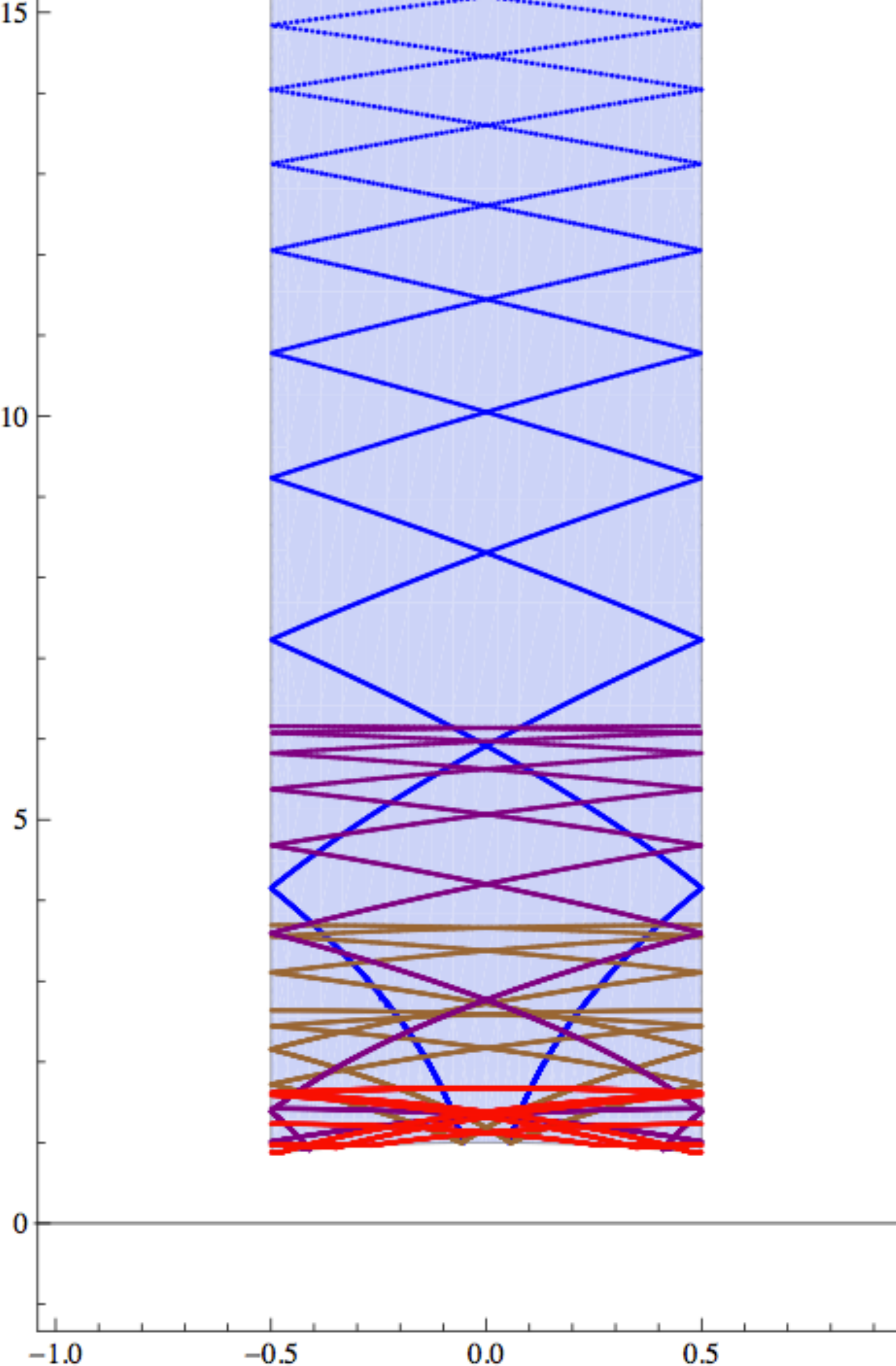}
                \caption{$D=1365:$ All} 
                \label{AKfig:D1365}
        \end{subfigure}
\caption{The  geodesics in $\sC_{D}$ corresponding to \eqref{AKeq:D1337gas} and \eqref{AKeq:D1365gas}.} 
\label{AKfig:geosDuke}
\end{figure}

For $D=
1337$ (which is $\equiv1\pmod4$ and square-free, and hence fundamental), we find that the class number is $h_{1337}=2$ and the class group $\sC_{1337}$ is comprised of two classes, represented by $Q_{1}=[7,-35,-4]$ and $Q_{2}=[4,-35,-7]$. 
Following the duality above to hyperbolic matrices, we find that $Q_{1}$ and $Q_{2}$ correspond, respectively, to
$$
M_{1}
\ = \
\left(
\begin{array}{cc}
 2\ 676\  336\  167 & 
 523\  561\  808 
 \\
 299\  178\  176 
 & 58\  527\  127 \\
\end{array}
\right)
,
\quad\text{and}\quad
M_{2}
\ = \
M_{1}^{\top}
.
$$
The entries are so large because the fundamental solution $(T,S)=(t,s)$ to $T^{2}-1337S^{2}=4$ is:
$$
(t,s) \ = \ (2\  734\  863\  294,\ 74\  794\  544),
$$
which also explains why the class group is so relatively small (cf. the discussion below \eqref{AKeq:hDis}).
The larger roots $\ga_{1}$ and $\ga_{2}$ of $Q_{1}(x,1)$ and $Q_{2}(x,1)$, respectively, have continued fraction expansions:
\be\label{AKeq:D1337gas}
\ga_{1} \ = \ [\overline{8, 1, 17, 2, 1, 1, 3, 1, 35, 1, 3, 1, 1, 2, 17, 1, 8, 5}]
,
\ee
$$
\ga_{2} \ = \ [\overline{5,8, 1, 17, 2, 1, 1, 3, 1, 35, 1, 3, 1, 1, 2, 17, 1, 8}]
.
$$
These have the same even length, and  differ by 
an odd
shift (by one), so the two forms $Q_{1}$ and $Q_{2}$ {\it are} equivalent under $\GL_{2}(\Z)$-action but not under $\SL_{2}(\Z)$. That is, while the narrow class number of 1337 is $2$, the full class number is $1$.
The corresponding geodesics are actually the same curve, but with opposite orientation; so they
look
 the same in $\G\bk\bH$ but differ in $T^{1}(\G\bk\bH)$, because the tangent vectors all change direction.
See  \figref{AKfig:D1337} for an illustration, which appears as just one curve, since the tangent vectors are not displayed.

For another example, we study the case $D=1365$. This discriminant is also $\equiv1\pmod4$ and square-free, and therefore fundamental. This time the class number is $h_{1365}=8$ ({\bf Exercise}), and we also leave it an exercise to work out representatives for the classes $[Q_{1}],\dots,[Q_{8}]$, the corresponding hyperbolic matrices $M_{1},\dots, M_{8}$, and the Pellian solution $(t,s)$. 
The roots $\ga_{1},\dots,\ga_{8}$ have continued fractions:
\be\label{AKeq:D1365gas}
\ga_{1}\ = \ [\overline{1, 35}],\
\ga_{2}\ = \ [\overline{35, 1}],\
\ga_{3}\ = \ [\overline {5, 7}], \
\ga_{4}\ = \ [\overline{7, 5}], 
\ee
$$
\ga_{5}=[\overline {1, 1, 1, 11}], \
\ga_{6}=[\overline {1, 1, 11, 1}], 
$$
$$
\ga_{7}=[\overline {1, 1, 1, 2, 1, 2}], \
 \ga_{8}=[\overline{1, 1, 2, 1, 2, 1}].
$$
Each pair $\ga_{2j-1},\ga_{2j}$ is the same geodesic but with opposite orientation;
the four distinct geodesics in $\G\bk\bH$ are illustrated in \figsref{AKfig:D136512}--\ref{AKfig:D136578}. Recall from \eqref{AKeq:TSto12s} that the cutting sequence of the geodesic flow is a symbolic coding of the visual point; then the first pair of $\ga$'s above corresponds to a geodesic that simply shoots up in the air and falls back down (\figref{AKfig:D136512}), the next two pairs are mid-level geodesics (\figsref{AKfig:D136534} and \ref{AKfig:D136556}), and the last pair of $\ga$'s gives a very low-lying geodesic (\figref{AKfig:D136578}). 

So the behavior of each geodesic in $\sC_{1365}$ is quite different, while there is only 
one geodesic (up to orientation) in $\sC_{1337}$. That said, 
combining all four geodesics in $\sC_{1365}$ into one picture (\figref{AKfig:D1365}),
%
one finds that, were the geodesics in $\sC_{1365}$ all colored the same, it would be quite difficult to distinguish
 this image from \figref{AKfig:D1337} for $\sC_{1337}$. Moreover the density of both plots is reminiscent of the invariant measure $dx\, dy/y^{2}$ on $\bH$ from \eqref{AKeq:Haar}. 
That is, these curves, when grouped by discriminant, become ``equidistributed'' 
with respect to the invariant measure
as the discriminant grows.
Equidistribution here means that the average amount of time this union of curves spends in a given nice region, say $A\subset\G\bk\bH$, becomes proportional to the area of the region.
This observation 
was
turned into a beautiful  theorem by Duke in 
\cite{Duke1988}. To formulate the statement more precisely, suppose such an $A$ is given, and write $\bo_{A}$ for the indicator function of $A$ in $\G\bk\bH$.
\begin{theorem}[Duke's Theorem\index{Duke's Theorem}]
As $D\to+\infty$ through fundamental discriminants,
\be\label{AKeq:DukeThm}
\frac1{h_{D}}
\sum_{\g\in \sC_{D}}
\frac1{\ell(\g)}
\int_{\g}\bo_{A}\,ds
\ \longrightarrow \
\frac1{\vol(\G\bk\bH)}
\int_{\G\bk\bH}\bo_{A}\,{dx\, dy\over y^{2}}.
\ee
Here $\ell(\g)$ is the length of the closed geodesic $\g$, and $ds$ is (hyperbolic) arclength measure.%
\end{theorem}
For simplicity, 
we have stated \eqref{AKeq:DukeThm} for the base space $\G\bk\bH$, but a similar result holds for the unit tangent bundle as well. Also, one can prove \eqref{AKeq:DukeThm} with effective (power savings) rates, and dropping the ``fundamental'' condition, see, e.g., \cite{ClozelUllmo2004}.
Again, this is a big theory of which we have only scratched the surface;
the proceeding discussion will suffice for our purposes here.

\section{Lecture 2: Three Problems in Continued Fractions: ELMV, McMullen, and
  Zaremba}


\subsection{
ELMV}

Duke's proof  of \eqref{AKeq:DukeThm} is a tour-de-force of analytic gymnastics, involving
a 
Maass
``theta correspondence''
to convert period integrals into Fourier coefficients of half-integral weight forms, and implementing methods pioneered by Iwaniec \cite{Iwaniec1987} (using Kuznestov's formula, sums of Kloosterman sums, and estimates of Bessel-type functions) to give non-trivial estimates of such. 
Meanwhile, the theorem itself \eqref{AKeq:DukeThm} seems like a beautifully simple dynamical statement; 
perhaps there is
a more ``ergodic-theoretic'' proof? Indeed, one was eventually 
obtained, after much work,
by 
Einsiedler-Lindenstrauss-Michel-Venkatesh 
 \cite{ELMV2012}
 (henceforth referred to as ELMV), with 
 many
 ideas
preempted decades earlier by Linnik \cite{Linnik1968}. 
 
ELMV wanted to approach
the problem (and higher rank analogs \cite{ELMV2009}) from a type of ``measure rigidity'' \'a la Rather's theorems -- one must show that 
the 
only
measure arising as a limit of the measures on the left hand side of \eqref{AKeq:DukeThm} is the 
Haar
measure $dx\, dy/y^{2}$. This raised the question: does one really need the full {\it average} over the class group here, or could individual geodesics already equidistribute?
Despite some partial progress to the affirmative (see, e.g., \cite{HarcosMichel2006, Popa2006}),
because of the symbolic coding of the geodesic flow,
 closed 
geodesics can be made to
have arbitrary behavior, simply by 
 choosing the partial quotients
 in the visual point and working backwards. So certainly not all long closed geodesics equidistribute. But perhaps  if one restricts only to {\it fundamental} closed geodesics, that is, ones whose corresponding discriminant is fundamental, the equidistribution is restored? 
Even then it is easy to produce examples of   non-equidistributing sequences of closed geodesics,
which, e.g., have limit measure $dy/y$ supported on the imaginary axis (see \cite[p. 233]{Sarnak2007}).
But what if we do not allow the ``mass'' to escape into the cusp of $\G\bk\bH$? Can we find closed fundamental geodesics which stay away from the cusp? Among these would certainly be some interesting limiting measures, and would demonstrate the difficulty of Duke's Theorem.
Such considerations naturally led ELMV around 2004 to propose the following 

\

{\bf Problem: (ELMV)}
Does there exist a compact subset $\cY\subset\cX=T^{1}(\G\bk\bH)$ of the unit tangent bundle of the modular surface which contains infinitely many fundamental closed geodesics?

\

The answer to this question turns out to be YES, as resolved (in a nearly best-possible quantitative sense)  by Bourgain and the author in \cite{BourgainKontorovich2015}. We will say more about the proof in the last lecture, but first move on to two other (seemingly unrelated) problems.

\subsection{McMullen's (Classical) Arithmetic Chaos Conjecture\index{Arithmetic Chaos Conjecture}}

Here is
a more recent problem posed by McMullen \cite{McMullen2009, McMullenNotes}, which he calls ``Arithmetic Chaos'' 
(not to be confused with Arithmetic Quantum Chaos, for which see, e.g. 
\cite{Sarnak2011a}.) 
The problem's statement begins in a similar way to ELMV, asking for closed geodesics on $\cX$ contained in a fixed compact set $\cY$, 
but the source is quite different.

McMullen was studying questions around the theme of
 Margulis's conjectures on the rigidity of higher-rank torus actions, 
and observed
 that
a very interesting problem in rank 1 
had been overlooked.\footnote{%
See \cite[Conj. 6.1]{McMullen2009} for a more precise statement which implies Arithmetic Chaos, and moreover predicts that the entropy discussed below can be made arbitrarily close to the natural limit.
} 
The statement is the following.

\begin{conjecture}[Arithmetic Chaos]\label{AKconj:ACI}
There is a compact subset $\cY$ of $\cX$ such that, for all real quadratic fields $K$, the set of closed 
geodesics defined over $K$ and lying in $\cY$ has positive entropy.
\end{conjecture}

Before explaining the meaning of the above words, we reformulate the conjecture as a simple statement about continued fractions.

\begin{conjecture}[Arithmetic Chaos II]\label{AKconj:ACII}
There is an $A<\infty$ so that, for any real quadratic field $K$, the set
\be\label{AKeq:ACII}
\{[\overline{a_{0},a_{1},\dots, a_{\ell}}]\in K : \text{ all }a_{j}\le A\}
\ee
has
exponential growth (as $\ell\to\infty$).
\end{conjecture}

The reformulation is quite simple. Recall yet again that  the geodesic flow is a symbolic coding of the  continued fraction expansion of the visual point. Thus going high in the cusp means having large partial quotients, and vice-versa. So a geodesic which is ``low-lying'' in some compact set $\cY$, that is, avoiding the cusp, is one whose visual point has only small partial quotients. 
Now since closed geodesics correspond to classes of binary forms with roots that are quadratic irrationals, the visual points automatically lie in some real quadratic field; in fact, it is easy to see that they lie in $K=\Q(\sqrt D)$, where $D$ is the discriminant of the geodesic (i.e. that of the class). This explains the appearance of real quadratic fields in both versions of the conjecture, as well as how being ``low-lying'' in \conjref{AKconj:ACI} corresponds to having small partial quotients in \conjref{AKconj:ACII}.
Without defining entropy, let us simply say that this condition in the first version corresponds in the second to the exponential growth of the set in \eqref{AKeq:ACII}.

In the third lecture below, we will present a certain ``Local-Global Conjecture''  which would easily imply McMullen's in the strongest from, that is, with $A=2$ in \conjref{AKconj:ACII}.
The same conjecture also has as an immediate consequence the aforementioned resolution of the ELMV Problem, as well as Zaremba's Conjecture (to be described below).
Unlike the latter two problems, it seems McMullen's problem requires the full force of this Local-Global Conjecture;
embarrassingly, the only meager progress made so far is numerical, as we now describe.

In McMullen's lecture \cite{McMullenNotes}, he gives numerical evidence for \conjref{AKconj:ACII} with $A=2$, taking the case $K=\Q(\sqrt5)$: he is able to find the following (primitive, modulo cyclic permutations and reversing of partial quotients) continued fractions:
$$
[\overline{1}]\ = \ \frac{1+\sqrt{5}}{2},\qquad [\overline{1, 1, 1, 1, 1, 1, 2, 1, 1, 2, 2, 1, 1, 1, 1, 2, 2}] \ = \ 
\frac{554+421 \sqrt{5}}{923}.
$$
McMullen presents these two surds as 
evidence of 
``exponential'' growth.

Using the Local-Global Conjecture as a guide, the author
 found the following further examples:
$$
[\overline{1, 1, 1, 2, 1, 2, 2, 2, 2, 1, 1, 2, 2, 1, 2, 1, 1, 2, 2, 1}] \   =\ \frac{90603+105937 \sqrt{5}}{207538} 
,
$$
$$
[\overline
{2, 1, 1, 2, 1, 1, 1, 1, 2, 2, 1, 1, 1, 1, 1, 1, 1, 1, 1, 2}
]
\ = \ 
\frac{12824 + 7728 \sqrt5 }{11667}
.
$$
Most recently, Laurent Bartholdi and Dylan Thurston (private communication) have pushed these numerics even further, finding the following further examples:
$$
[\overline
{1, 1, 
1, 1, 1, 1, 2, 2, 2, 1, 1, 1, 1, 1, 2, 1, 1, 1, 2, 2, 2, 1, 2, 2, 1, 1, 1, 2
}
],
$$
$$
[\overline
{1, 1, 1, 1, 1, 2, 2, 1, 1, 1, 2, 1, 
1, 1, 1, 1, 2,2, 2, 1, 2, 1, 1, 1, 1, 2, 2, 2}],
$$
$$
[\overline
{1, 1, 1, 2, 1, 1, 2, 1, 2, 2, 2, 2, 2, 2, 1, 2, 1, 1, 2, 1, 2, 1, 1, 2, 
1, 2, 2, 2, 2, 2, 2, 1, 2, 2}
],
$$
$$
[\overline
{1, 1, 1, 1, 1, 2, 2, 1, 1, 1, 2, 2, 1, 2, 1, 2, 2, 1, 2, 2, 2, 1, 1, 2, 
2, 2, 2, 2, 1, 2, 2, 1, 2, 1, 2, 2}],
$$
$$
[\overline
{1, 1, 1, 1, 2,1, 1, 2, 2, 1, 1, 1, 2, 2, 1, 1, 1, 2, 2, 2, 1, 2, 1, 1, 2, 2, 2, 1, 
1, 1, 2, 1, 1, 1, 2, 2}],
$$
$$
[\overline
{1, 1, 2, 1, 2, 1, 1, 2, 2, 2, 1, 1, 2, 1, 2, 2, 2, 2, 1, 1, 2, 2, 1, 2, 
1, 1, 2, 2, 2, 1, 1, 2, 1, 2, 2, 2}],
$$
$$
[\overline
{1, 1, 1, 1, 1, 1, 1, 1, 2, 1, 1, 1, 1, 1, 2, 2, 2, 2, 2, 2, 1, 2, 2, 1, 
2, 2, 2, 1, 1, 2, 1, 2, 1, 2, 2, 2, 2, 1, 2, 1, 1, 2}
],
$$
$$
[\overline
{1, 1, 1, 1, 1, 1, 1, 1, 1, 2, 1, 1, 2, 1, 1, 1, 1, 1, 2, 1, 2, 2, 2, 2, 
2, 1, 1, 2, 1, 2, 2, 2, 2, 2, 1, 2, 1, 2, 2, 1, 2, 2, 2}],
$$
$$
[\overline
{1, 1, 1, 2, 1, 1, 2, 1, 2, 2, 1, 1, 1, 2, 1, 2, 1, 1, 2, 1, 1, 1, 2, 2, 
1, 1, 1, 2, 2, 2, 2, 2, 1, 1, 1, 2, 2, 2, 2, 1, 1, 1, 2, 2}].
$$
all in $\Q(\sqrt 5)$.
One may now argue whether this list (of 13 distinct surds in all) is yet indicative of exponential growth.

\subsection{Zaremba's Conjecture\index{Zaremba's Conjecture}}

Our final problem originates in questions about pseudorandom numbers and numerical integration. A detailed discussion of these questions is given in \cite[\S2]{Kontorovich2013}, so we will not repeat it here. The statement of the conjecture, understandable by Euclid, is as follows.\footnote{%
See \cite[\S6]{McMullen2009} for McMullen's connection of Zaremba's Conjecture to Arithmetic Chaos.
}
\begin{conjecture}[Zaremba \cite{Zaremba1972}]
There is some $A<\infty$ such that, for every integer $d\ge1$, there is a coprime integer $b\in(0,d)$ such that
the reduced rational $b/d$ has the (finite) continued fraction expansion
$$
\frac bd = [0,a_{1},\dots, a_{\ell}]
,
$$
with all partial quotients $a_{j}$ bounded by $A$.
\end{conjecture}

Progress on the advertised Local-Global Conjecture allowed Bourgain and the author to resolve a 
density
version of Zaremba's Conjecture:

\begin{theorem}[\cite{BourgainKontorovich2014}]\label{AKthm:BKZarem}
There is an $A<\infty$ such that the proportion of $d<N$ for which Zaremba's conjecture holds approaches $1$ as $N\to\infty$. 
\end{theorem}

In the original paper, $A=50$ was sufficient, and this has since been reduced to $A=5$ in \cite{Huang2015, FrolenkovKan2014}. Most recently, Zaremba's conjecture has found application to counterexamples to Lusztig's conjecture on modular representations, via the groundbreaking work of Geordie Williamson, see \cite{Williamson2015}.

\section{Lecture 3: The Thin Orbits Perspective}

All three problems discussed in Lecture 2 are collected here under a common umbrella as a
``Local-Global'' problem for certain  ``thin'' orbits. We discuss recent joint
work with Jean Bourgain which settles the first problem (ELMV), and
makes some progress towards the last (Zaremba).

To motivate the discussion, recall from the Exercise below \eqref{AKeq:MfromAs} and the expression \eqref{AKeq:gaIs} that the quadratic surd 
$$
\ga=[\overline{a_{0},a_{1},\dots, a_{\ell}}]
$$
is fixed by the matrix
$$
M=
\mattwo {a_{0}}110\cdot
\mattwo {a_{1}}110\cdots
\mattwo {a_{\ell}}110
,
$$
with
$$
\ga\ \in \ \Q(\sqrt{\tr^{2}M-4}).
$$
This tells us that studying traces of such matrices $M$, we might learn something about both ELMV and McMullen's conjectures. For example:

{\bf Exercise:} If $\tr^{2}M-4$ is squarefree, then the corresponding closed geodesic is fundamental.
 
Yet another elementary exercise is that, if $b/d=[0,a_{1},\dots, a_{\ell}]$, then
$$
\mattwo {a_{1}}110\cdots
\mattwo {a_{\ell}}110
=
\mattwo db**
,
$$ 
so studying the 
top-left entries of matrices of the above form tells us about Zaremba's conjecture.

In all the above problems, the partial quotients $a_{j}$ are bounded by some absolute constant $A$. Thus we should study the semigroup generated by matrices of the form $\mattwos a110$, with $a\le A$. 
Actually, we will want all elements in $\SL_{2}(\Z)$ (whereas the generators have determinant $-1$), 
 so we define the key semigroup
$$
\G_{A}\  : = \ \<\mattwo a110 \ : \ a\le A\>^{+} \ \cap \ \SL_{2},
$$
of
even length words in the generators; here
the superscript $+$ denotes generation as a semigroup. 

We first claim that this semigroup $\G_{A}$ is {\it thin}. We give the definition  by example. As soon as $A\ge2$, the Zariski closure of $\G_{A}$ is $\SL_{2}$, that is, the zero set of all polynomials vanishing on all of $\G_{A}$, is 
that of the single polynomial $P(a,b,c,d)=ad-bc-1$.
The  integer points, $\SL_{2}(\Z)$, of the Zariski closure grow like:
$$
\#(\SL_{2}(\Z)\ \cap \ B_{X}) \ \asymp \ X^{2},
$$
where $B_{X}$ is a ball about the origin of radius $X$ in one's favorite fixed archimedean norm. 
On the other hand, an old result of Hensley \cite{Hensley1989} 
gives that
\be\label{AKeq:GamAsize}
\#(\G_{A} \ \cap\ B_{X}) \ \asymp \ X^{2\gd_{A}}.
\ee
Here $\gd_{A}$ is the Hausdorff dimension of the ``limit set'' $\fC_{A}$ of $\G_{A}$, defined as follows:
$$
\fC_{A} \ := \ \{[0;a_{1},a_{2},\dots]\ : \ \forall j, \ a_{j}\le A\}.
$$
The dimensions of these Cantor-like sets have been studied for a long time \cite{Good1941, JenkinsonPollicott2001}, e.g.,
\be\label{AKeq:del2Is}
\gd_{2} \  \approx  \ 0.5312, \qquad
\gd_{3} \  \approx  \ 0.7056, \qquad
\gd_{4} \  \approx  \ 0.7889, 
\ee
and as $A\to\infty$, Hensley \cite{Hensley1992} showed that
$$
\gd_{A} \ = \ 1-{6\over \pi^{2}A}+o\left(\frac1A\right).
$$
The point is that all of these dimensions $\gd_{A}$ are strictly less than $1$, so it follows from \eqref{AKeq:GamAsize} that
$$
\#(\G_{A}\cap B_{X}) \ = \ o\left(\#(\SL_{2}(\Z)\cap B_{X})\right),
$$
as $X\to\infty$. 
This is the defining feature of a so-called {\it thin integer set}\index{thin integer set} (defined this way in \cite[p. 954]{Kontorovich2014}) -- it has archimedean zero density in the integer points of its Zariski closure. When the set in question is actually a  sub{\it group} of a linear group, this definition agrees with the ``other'' definition of thinness, namely the infinitude of a corresponding co-volume, or index; see, e.g., \cite{Sarnak2014}.
\\

The ELMV, McMullen and Zaremba problems do not just study the semigroup $\G_{A}$ itself;
they all study the image of $\G_{A}$ under some 
%
linear map $F$ on $M_{2\times 2}(\Z)$ taking integer values on $\G_{A}$. 
For 
example, 
in
Zaremba's conjecture, one takes $F:\mattwos abcd\mapsto a$ which picks off the top-left entry. For ELMV and Arithmetic Chaos, one studies the trace, $F:\mattwos abcd\mapsto a+d$. Then 
the key object of interest is $F(\G_{A})\subset\Z$, and even more precisely, the multiplicity with which an integer $n$ is represented in $F(\G_{A})$. We define the multiplicity as:
$$
\mult (n) \ := \ \#\{\g\in\G_{A} \ : \ F(\g)=n\}. 
$$
{\it A priori} this count may be infinite (e.g., if $F$ is constant), so it is useful to also define a quantity guaranteed to be finite, by truncating:
$$
\mult_{X}(n) \ : = \ \#\{\g\in\G_{A}\cap B_{X} \ :  \ F(\g) = n\}.
$$
Since $F$ is 
linear, 
the image of an archimedean ball $B_{X}$ will also be of order $X$, so one may 
naively expect the multiplicity of some $n\asymp X$ to be of order 
\be\label{AKeq:heur}
\frac1X \cdot \#(\G_{A}\cap B_{X}) \ \asymp \ X^{2\gd_{A}-1}.
\ee
It is easy to see that such a prediction is too primitive, e.g., if $F:\mattwos abcd\mapsto 2a$, then all odd integers are missing in the image. 

Given $\G_{A}$ and $F$, we define an integer $n\in\Z$ to be {\it admissible} if it passes all ``local obstructions'':
$$
n \ \in \ F(\G_{A})\pmod q,
$$
for every integer $q\ge1$. While this  condition may seem difficult to verify (for instance, it asks about infinitely many congruences!), it turns out that, thanks to the theory of Strong Approximation, it is very easy to check in practice; see, e.g., \cite[\S2.2]{Kontorovich2013} for a discussion.

The Local-Global Conjecture\index{Local-Global Conjecture}, formulated by
Bourgain and the author  states that, once these local obstructions are passed, the naive heuristic \eqref{AKeq:heur} holds.

\begin{conjecture}[The Local-Global Conjecture]
Assume that $A\ge2$, so that $\G_{A}$ is Zariski dense in $\SL_{2}$, and that its image under the 
linear map  $F$ is infinite; equivalently the Zariski closure of $F(\G_{A})$ is the affine line. 
For a growing parameter $X$ and an integer $n\asymp X$ which is admissible, we have
\be\label{AKeq:LocGlob}
\mult_{X}(n) \ = \ X^{2\gd_{A}-1-o(1)}.
\ee
\end{conjecture}

Notice that the dimensions in \eqref{AKeq:del2Is} all exceed $1/2$, whence the exponents $2\gd_{A}-1$ in \eqref{AKeq:LocGlob} are all positive. So for large $X$, these multiplicities are non-zero, that is, large numbers which pass all local obstructions should be ``globally'' represented in $F(\G_{A})$. 

We leave it as a pleasant {\bf Exercise} to prove McMullen's \conjref{AKconj:ACII} from \eqref{AKeq:LocGlob}; see also \cite[Lemma 1.16]{BourgainKontorovich2013a}.

The progress leading to \thmref{AKthm:BKZarem}  uses the fact that the Zaremba map $F:\mattwos abcd\mapsto a$ is of ``bilinear type,'' in that $F$ can be written as $F(M)=\<v_{1}\cdot M,v_{2}\>$, where
$v_{1}=v_{2}=(1,0)$. 
A similar result can be proved (see \cite{BourgainKontorovich2010}) whenever $F$ is of this form; equivalently ({\bf Exercise}) whenever $\det F = \ga \gd - \gb \g =0$, where
$F:\mattwos abcd\mapsto \ga a + \gb b +\g c + \gd d$.


\begin{figure}
\includegraphics[width=\textwidth]{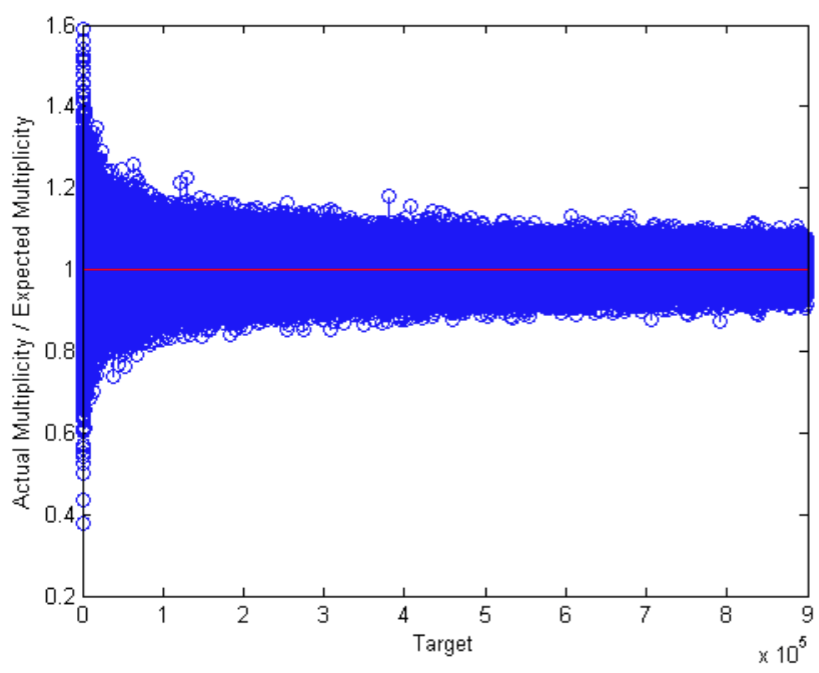}
\caption{Numerical verification of  \conjref{AKconj:Zar} for $A=5$. Image by P. Cohen in \cite{Cohen2015}.}
\label{AKfig:Zaremb}
\end{figure}

Quite recently, a  more precise formulation of Zaremba's conjecture, based on the Local-Global Conjecture, and the execution of certain Hardy-Littlewood-type estimations and identities for the ``singular series'' was given by Cohen \cite{Cohen2015}, a Rutgers summer REU student. 
\begin{conjecture}\label{AKconj:Zar}
In the Zaremba setting of $F:\mattwos abcd\mapsto a$, we have
$$
\mult(n) \ \sim \ 2\gd_{A}{\#(\G_{A}\cap B_{n})\over n}\times {\pi^{2}\over 6}\times \prod_{p\mid n}\left(1-{1\over p}\right).
$$
\end{conjecture}

A plot of the left hand side divided by the right hand side is given in \figref{AKfig:Zaremb}; the data is rather convincing in support of the refined conjecture, which hopefully also serves as evidence for the more general Local-Global Conjecture.

\subsection{Tools: Expansion and Beyond}\

We give here but a hint of some of the methods developed to prove the results in \cite{BourgainKontorovich2014, BourgainKontorovich2015}.
A key initial ingredient is what we shall refer to broadly as ``expansion,'' or ``SuperApproximation.'' This has also been discussed rather extensively in numerous surveys, see, e.g., \cite{Sarnak2004, Lubotzky2012}, and \cite[\S5.2]{Kontorovich2013}, \cite[\S3.3]{Kontorovich2014}. To go ``Beyond Expansion,'' one needs to develop more sophisticated tools outside the scope of this survey; a hint is given in \cite[\S3.7]{Kontorovich2014}, and the interested reader is invited to peruse the original papers \cite{BourgainKontorovich2015a, BourgainKontorovich2015, BourgainKontorovich2016}.
\\

\textbf{Acknowledgements}
 The author is grateful to the organizers for the opportunity to collect these topics, and to Jean Bourgain for the collaborative work described here. Thanks also to Valentin Blomer and Curt McMullen for comments on an earlier draft.

\bibliographystyle{alpha}

\bibliography{../../AKbibliog}

\end{document}